\documentclass[final,3p,times]{elsarticle}

\usepackage{amssymb}

\usepackage{amsmath}
\usepackage{booktabs}
\usepackage{caption}
\usepackage{subcaption}
\usepackage{array}
\usepackage{xcolor}
\usepackage{textcomp}
\usepackage{soul}

\biboptions{comma, square, sort&compress}

\journal{Applied Mathematical Modelling}

\begin{document}

\begin{frontmatter}



\title{Discrete Transparent Boundary Conditions\\for the Equation of Rod Transverse Vibrations}


\author{Vladimir A. Gordin}
\ead{vagordin@mail.ru}

\author{Aleksandr A. Shemendyuk\corref{shemfn}}
\cortext[shemfn]{Corresponding author}
\ead{aashemendyuk@edu.hse.ru}

\address{National Research University "Higher School of Economics", Moscow, Russia}
\address{Hydrometeorological Centre of Russia, Moscow, 123242, Russia}

\begin{abstract}
Local perturbations of an infinitely long rod travel to infinity. On the contrary, in the case of a finite length of the rod, the perturbations reach its boundary and are reflected. The boundary conditions constructed here for the implicit difference scheme imitate the Cauchy problem and provide almost no reflection. These boundary conditions are non-local with respect to time, and their practical implementation requires additional calculations at every time step. To minimise them, a special rational approximation, similar to the Hermite - Pad\'e approximation is used. Numerical experiments confirm the high ``transparency'' of these boundary conditions {and determine the conditional stability regions for finite-difference scheme}.
\end{abstract}


\end{frontmatter}

\section{Introduction}
\label{S:1}

The equation of transverse vibrations of a rod (beam) with a circular cross section 
\begin{equation}
\label{1.1}
\rho\frac{\partial^2 u}{\partial t^2}-\frac{\partial}{\partial x}\left[R^2\rho\frac{\partial^3 u}{\partial x \,\partial t^2}\right] +\frac{\partial^2 }{\partial x^2}\left[ER^2 \frac{\partial^2 u}{\partial x^2}\right]=f
\end{equation}
has many applications, see e.g. \cite{timoshenko2008vibration}. Here {$x\in [-L/2,\,L/2]$ is the spatial variable, $L$ is the length of the rod, $t$ is time}, $\rho$ is the density of the rod, $R$ is the radius of the cross section of the rod, $E$ is the Young's modulus of the rod, the unknown function $u=u(t,\,x)$ describes the transverse displacement of the rod. The right hand side (forcing) $f$ describes external force.

{The kinetic energy of transverse vibrations of a rod is defined as}
\begin{equation}
    \label{eq:kinetic_energy}
    \mathbf{K}[u] = \frac{1}{2} \int\limits_{-L/2}^{L/2}\rho\left[\left(\frac{\partial u}{\partial t} \right)^2+R^2\left(\frac{\partial^2 u}{\partial t \, \partial x} \right)^2\right] \,\mathrm{d}x,
\end{equation}
{and its potential energy is}
\begin{equation}
    \label{eq:potential_energy}
    \mathbf{P}[u] = \frac{1}{2}\int\limits_{-L/2}^{L/2}ER^2\left(\frac{\partial^2 u}{\partial x^2} \right)^2\,\mathrm{d}x.
\end{equation}

{Eq.~\eqref{1.1} could be obtained from the least action variation principle with Lagrangian $\mathcal{L} = \mathbf{K} - \mathbf{P}$. The Hamiltonian (energy) for Eq.~\eqref{1.1} is then} 
\begin{equation}
    \label{norm:H}
    \mathcal{H} = \mathbf{K} + \mathbf{P}.
\end{equation}

The Cauchy problem for Eq.~(\ref{1.1}):
\begin{equation*}
u(0,x)=U_0(x),\quad\partial_t u(0,x)=U_1(x),\quad x\in\mathbb{R},
\end{equation*}
is correct according to Hadamard. Two boundary conditions at each edge of the segment $x=\pm L/2$ provide the mixed initial-boundary value problem, if the corresponding boundary operators satisfy the Shapiro -- Lopatinsky conditions {(see~\cite{agranovich1971boundary, hormander2007analysis, Kreiss1970, KreissHedwig2001, sakamoto1970i, sakamoto1970ii})}.

{Let us assume that the forcing in Eq.~\eqref{1.1} is absent: $f \equiv 0$.} The energy $\mathcal{H}$ remains unchanged for the Cauchy problem for Eq.~\eqref{1.1}. If, in the Cauchy problem, the support of initial value functions belongs to the segment
\begin{equation*}
    \mathrm{supp} \, U_0 \subseteq [-L/2, \, L/2], \quad \mathrm{supp} \, U_1 \subseteq [-L/2, \, L/2],
\end{equation*}
{then for any time $t$ the following inequality holds:}
\begin{equation}
    {\bf P}[u(t, x)]\le {\cal H}[u(t,\,x)]\le \frac{1}{2} \int\limits_{-\infty}^{\infty} \left[\rho \left(U_1\right)^2 + \rho  R^2 \left(\frac{\partial  U_1}{ \partial x} \right)^2+ ER^2\left(\frac{\partial^2 U_0}{\partial x^2} \right)^2\right] \, \mathrm{d}x,
\end{equation}
{because the sub-integral expression is always non-negative. In numerical approximations of Eq.~\eqref{1.1} the first integral $\mathcal{H}$ can change with time due to rough discretization or because of boundary conditions.}

{
The construction of boundary conditions that have no reflection of outgoing waves from the boundary for an implicit finite-difference equation (see Sect.~\ref{se.schem}), which approximates Eq.~(\ref{1.1}), is the subject of this article. There are various names of such boundary conditions (see Sect.~\ref{S:12}). We use the name DTBCs (Discrete Transparent Boundary Conditions). We emphasise that for each finite-difference scheme it is necessary to construct ``individual'' DTBCs.} 

{
The construction is practically important, because of the place of this equation in the theory of elasticity and in engineering problems. The algorithm of the construction is technically more consuming than in the case of classical equations and systems of mathematical physics in partial derivatives considered earlier in {\cite{engquist1977absorbing, gordin1977, engquist1979radiation, gordin1978mixed, gordin1979diss, gordin1987mathematical, gordin2000mathematical, arnold2005approximation, gordin2010mathematics, arnold2003}}, since \textbf{a)} Eq.~(\ref{1.1}) is not resolved with respect to the highest (second) derivative with respect to time (i.e. it is not a differential equation of the Cauchy -- Kovalevskaya type) and \textbf{b)} the equation's order with respect to the spatial variable is equal to $4$.}

{
The high order of the differential (and, as a consequence, of the finite-difference) equation with respect to space requires two boundary conditions at each edge. The DTBCs {and Approximate DTBC (or simply ADTBCs)} for various simpler finite-difference schemes approximating equations of mathematical physics were considered {by authors} earlier in \cite{gordin2000mathematical, gordin2010mathematics, gordin1978projectors, gordin1979upper, gordin1982boundary, gordin1982application, gordin1987mathematicalb}.}

{
The various versions of the TBCs, DTBCs, and {ADTBCs}  may be developed for various equations and systems of the elasticity theory and their various discretizations, and such {ADTBCs} give the possibility of avoiding significant errors as a result of a false reflection from the boundary, when we solve local problems in a large area.}

{
We introduce the absolutely stable finite-difference scheme that approximates the Cauchy problem for Eq.~\eqref{1.1} with constant coefficients in Sect.~\ref{se.schem}. The construction of the DTBCs and {ADTBCs} for the finite-difference scheme will be considered in Sect.~\ref{se.TBC}.}

On the contrary, the implicit finite-difference scheme that approximates the mixed boundary value problem for Eq.~\eqref{1.1} is not absolutely stable {with our proposed ADTBCs}. Stability conditions significantly differ from similar conditions for schemes that approximate the classical equations of mathematical physics (transport, diffusion, wave, Schr\"odinger, etc.). {For them, usually, the stability condition is that} the time step $\tau$ must not exceed some critical value. According to our numerical experiments, for the considered mixed problem with given boundary conditions the step $\tau$ must satisfy two inequalities at once. {It appears} that on the $\tau$-$h$ plane the range of values, at which stability takes place lies between two parabolas, see Subsect.~\ref{subsect.stability}.

The description of the results of numerical experiments with {ADTBCs}, and comparisons with  other boundary conditions, are in Sect.~\ref{se.resul}.

\section{Transparent Boundary Condition Problem Overview}
\label{S:12}

In many practical problems of mathematical physics the necessary complete set of physically adequate boundary conditions is absent. For example, in meteorological models (see e.g. \cite{gordin1979diss, gordin1987mathematical, gordin2000mathematical}) we need  boundary conditions at the upper computational level $p=\mathrm{const} > 0$, but {no} boundary condition {can} describe {the variety of very complex meteorological phenomena that occur in the upper atmosphere (above the layer essential to the task of forecasting for one week). On the other hand, we cannot use the pressure level $p=0$ in the models, since we do not have steady measurements in the upper atmosphere, and the gas dynamics approximation becomes inadequate here. However, the upper boundary condition is necessary for the closure of both the differential and the finite-difference problems. Here the transparent boundary condition is a reasonable compromise. The corresponding {ADTBC} (see \cite{gordin1979upper}) decreases the forecasting error in comparison with ``simple'' boundary conditions.}

Also, in the problem of weather forecasting for a limited area $V$, we can use a ``simple'' approach and set Dirichlet boundary conditions at the border $\partial V$, where the right-hand side of the conditions is taken from a larger scale (global) forecasting model. 
We can consider the difference between meteorological fields in these two numerical models. The dynamics of this difference {(which is a vector-function)} show that waves coming out of the computational area are reflected from the boundary $\partial V$. It is not a physical effect, which worsens the regional model's forecasts.

To obtain a mathematically correct mixed initial-boundary problem for the differential system, the number $M$ of boundary conditions (according to the Shapiro –- Lopatinsky theory, see~\cite{agranovich1971boundary, Kreiss1970, KreissHedwig2001, sakamoto1970i, sakamoto1970ii, hormander2007analysis}) {may} be smaller than the number of unknown functions. Moreover, for the perfect gas (with constant heat capacity and without dissipation) dynamics, the number $M$ depends on the orientation of the wind direction at any boundary point $\vec x \in \partial V$. The number of boundary conditions for the corresponding dissipative (viscous) models is equal to the number of the unknown functions in the model.

\textbf{Note 1.} On the contrary, the Dirichlet boundary conditions for usual finite-difference approximations of the gas dynamic system provide the existence and uniqueness of the mixed problem's solution. However, the dissimilarity of these similar differential and finite-difference problems, significantly obstruct the convergence of finite-difference solutions when the steps of a finite-difference scheme tend to zero.

If we approximate a differential equation or a system using finite-difference schemes, the deficit of physically based boundary conditions may increase. For instance, if we use the finite-difference spatial approximation of the ideal gas dynamic's differential equations, according to the central difference formula, the number of boundary conditions required for the uniqueness of the solution increases in comparison with the differential problem. For the transport equation with viscosity, phenomena such as a boundary layer on the outflow from the region may occur.

Such computational difficulties are usually overcome by introducing into the algorithm a special finite-difference operator of a large computational viscosity in some vicinity of the boundary. However, it inevitably leads to an increase of the prediction error.

\textbf{Note 2.} To avoid the non-physical reflection in such boundary problems it is necessary to construct an analogue of the famous Sommerfeld radiation condition for the Helmholtz equation {$\Delta u + k^2u = 0$}, see e.g. \cite{sommerfeld1912greensche, Morawetz1965}. The asymptotical condition {(as $r \equiv |\vec x| \to \infty$)
\begin{equation*}
    \partial_r u - \mathrm{i} k u = \mathbf{o}\left(r^{-1}\right)
\end{equation*}
}
guarantees the solution's uniqueness of the Helmholtz equation in the whole space $\mathbb{R}^m$.

Some types of linear evolutionary partial differential equations may be reduced to stationary differential equations {(e.g., the wave equation to the Helmholtz one -- for the Cauchy's problem using the Laplace transform, or, if the solution is assumed to be harmonic with respect to time, using the Fourier transform)}. This asymptotical Sommerfeld condition may be developed for some other PDEs (and systems). 

The accuracy order of the asymptotical condition may be improved, if suitable terms with higher order derivatives are added. The conditions may be approximated by finite-difference formulae.

\textbf{Note 3.} The term Artificial (or Absorbing) Boundary Conditions (ABC) is also known as boundary conditions Imitating Cauchy Problem (ICP), {full absorption conditions}, computational boundary conditions, transparent boundary conditions, radiation boundary conditions, open boundary conditions, etc., see e.g. \cite{orlanski1976simple, bennett76}. The properties of such boundary conditions were discussed and compared in \cite{bayliss1980radiation, israeli1981approximation, keller1989exact, TSYNKOV1998465, hagstrom_1999}.

Let us consider a {linear differential equation with partial derivatives that is correct by Petrovsky (see \cite{Gelfand1967, gordin2010mathematics}):}
\begin{equation}
\label{gener}
\partial_t \vec u = A\vec u + \vec f,
\end{equation} 
where $A=A(\partial_{\vec x})$ is a linear differential operator in $\mathbb{R}^m$, $\vec f$ is a given function (forcing). We consider Eq.~(\ref{gener}) in area $V\subset \mathbb{R}^m$ with piece-smooth boundary $\partial V$. To determine a unique solution of  Eq.~(\ref{gener}) Cauchy initial conditions are posed:
\begin{equation*}
    \vec u(t,\,\vec x)|_{t=0} = \vec u_0(\vec x).
\end{equation*}

The given functions $\vec f\left(\vec x\right)$, $\vec u_0\left(\vec x\right)$ with supports in $V$, may be extended into the whole space by zeros and the Cauchy problem for Eq.~(\ref{gener}) may be considered for $\mathbb{R}^m$. Let us use this solution as a reference. {We say that a set of boundary conditions imitate the Cauchy problem, if for all $\vec u_0$ and $\vec f$ the solution of Eq.~(\ref{gener}) in the area $V$ under these boundary conditions is equivalent to the reference solution with extended by zero functions $\vec f\left(\vec x\right)$ and $\vec u_0\left(\vec x\right)$.}

\textbf{Note 4.} The definition may be extended for the cases of differential equations with a higher differential order with respect to time $t$ and, moreover, non-resolved with respect to the highest time derivative. Eq.~(\ref{1.1}) is an example of such a differential equation. {However, first we briefly describe TBCs for several classical linear PDEs {to facilitate understanding} (see e.g. \cite{gordin1977, gordin1978mixed, gordin1978projectors, gordin1982boundary, gordin1987mathematicalb, gordin2010mathematics, arnold2003, ryaben1990faithful, ryabenkii1993artificial, ryaben1995artificial, ryaben2006theoretical, fesch_2013}).}

{The reference solution may be obtained for such functions $\vec f,\:\vec u_0$, if the Green functions (fundamental solutions, kernels, etc.) of Eq.~(\ref{gener}) are known. Since ${\bf supp}\,\vec f\subseteq \bar V$ and ${\bf supp}\,\vec u_0\subseteq  \bar V$, we restrict the relative integrals over $\mathbb{R}^m$ in these formulae
on the area $V$ and obtain for $\vec x\in V$}:
\begin{equation}
\label{eqeq}
\vec u\left(t,\,\vec x\right) = \int\limits_V K_1 \left( t, \vec x-\vec y \right) \vec u_0 \left(\vec y\right)\, \mathrm{d}\vec y+ \int\limits_0^t\int\limits_V K_2 \left(t-s,\, \vec x-\vec y\right)\vec f\left(s, \vec y\right)\,\mathrm{d} \vec y \, \mathrm{d}s.    
\end{equation}

For instance, for the scalar diffusion equation $\partial_t u=D\Delta u + f$ the scalar  Green  functions are:
\begin{equation}
\label{kernconvol}
    K_1\left( t, \vec x-\vec y \right)=\frac{1}{\left(\sqrt{4\pi D t}\right)^m} \, \exp\left(\frac{- \left| \vec x-\vec y \right|^2}{4Dt}\right),\quad
    K_2 \left(t-s,\, \vec x-\vec y\right) = \frac{1}{\left(\sqrt{4\pi D(t-s)}\right)^m} \, \exp\left(\frac{-\left| \vec x-\vec y \right|^2}{4D(t-s)}\right).
\end{equation}

The famous Poisson integral formula, as well as integral formulae (convolutions) for other classical equations of mathematical physics with constant coefficients, may be obtained by the Fourier transform. They give an exact solution of the Cauchy problem for arbitrary initial data and right hand sides such that ${\bf supp}\,\vec f \subseteq V$, ${\bf supp}\,\vec u_0 \subseteq V$. {The transparent boundary conditions should provide the same solution for mixed problem in $x \in V, \; t \geq 0$.}

The kernels $K_1$ and $K_2$ decrease exponentially with respect to spatial variables. {However, the practical computation of integrals (\ref{eqeq}) may be computationally expensive.}

Although the Green functions exist for equations with variable coefficients, their explicit practical determination is very difficult. The finite-difference approach for the special mixed initial-boundary problem is more preferable, although its implementation is associated with noticeable difficulties.

{Sommerfeld's condition is fulfilled only  asymptotically as {$|\vec x|\to \infty$}. On the contrary, we define and construct boundary conditions {on $\partial V$, i.e. }in concrete points, lines, or planes. We consider the cases, where $\partial V$ is one or two points, one or two parallel lines, or planes.}

Let the variable $x_1$ be normal to the boundary, $\vec x=\langle x_1,\,\vec y\rangle,\;\mathrm{dim}\,\vec x=\mathrm{dim}\,\vec y+1$. To obtain TBC, we have to give up the locality property {(inherent, e.g., in the Sommerfeld condition for the wave equation)} and include non-local integral operators of the convolution type. Let $V=\mathbb{R}^m_+,\;x_1>0$. The integrals provide participation in such boundary conditions at an arbitrary time moment $t>0$ and in an arbitrary border point $\vec y \in \partial V$ of the solution's boundary values at all previous time moments $0 \leq s\le t$ and at $\vec z\in \partial V$, i.e. the boundary condition takes the form:
\begin{equation}
\label{kernconvol2}
    \partial_{x_1} u (t,\,0,\,\vec y) = \int\limits_{0}^t\int\limits_{\mathbb{R}^{m-1}}K\left(t-s,\,\vec y-\vec z \right)\,u\left(s,\,0,\,\vec z \right)\,\mathrm{d}\vec z\, \mathrm{d}s,
\end{equation}
{or
\begin{equation}
\label{kernconvol2-2}
    \partial_{x_1} u(t,\,0,\,\vec y) = Q(\partial_{\vec y}) u(t,\,0,\,\vec y) +
    \int\limits_{0}^t\int\limits_{\mathbb{R}^{m-1}} K\left(t-s,\,\vec y-\vec z \right)\,u\left(s,\,0,\,\vec z \right)\,\mathrm{d}\vec z\, \mathrm{d}s,
\end{equation}
or more general integro-differential operator
\begin{equation}
\label{kernconvol2-3}
\partial_{x_1} u(t,\,0,\,\vec y) = Q(\partial_{\vec y}) u(t,\,0,\,\vec y) + 
\sum\limits_{\kappa} \int\limits_{0}^t \int\limits_{\mathbb{R}^{m-1}} K_{\kappa}\left(t-s,\,\vec y-\vec z \right)\,
Q_{\kappa} u\left(s,\,0,\,\vec z \right)\,\mathrm{d}\vec z\, \mathrm{d}s,
\end{equation}
where $Q$ and $Q_{\kappa}$ are differential operators with respect to time and spatial variables that are tangential to the boundary. Here, we may have several terms with different differential operators and kernels -- the indexes of such terms are denoted by $\kappa$. The forms of these differential operators and kernels $K_{\kappa}$ are determined by the original operator $A$ (see examples, e.g., in \cite{ZAITSEV2007, gordin2010mathematics}).
}

{Boundary conditions Eq.~\eqref{kernconvol2}, \eqref{kernconvol2-2} and \eqref{kernconvol2-3} can imitate the Cauchy's problem with initial values and right hand side being prolonged with zeros. If the values of given functions $u_0(\vec x)$ and $f(t, \vec x)$ outside of the computational domain are nonzero, then the boundary condition Eq.~\eqref{kernconvol2}, \eqref{kernconvol2-2} and \eqref{kernconvol2-3} should be supplemented with terms that describe their contribution.}

The TBCs for the one dimensional wave equation
\begin{equation*}
    \partial_t^2 u =c^2\partial_x^2 u
\end{equation*}
on the segment $[-1,\,1]$ are local:
\begin{equation*}
    \partial_t u =\pm c \partial_x u.
\end{equation*}

Usually TBCs are non-local. For instance, for the diffusion equation in the half-space $x\ge 0,\,\vec y\in\mathbb{R}^{m-1}$ we obtain the following TBC:
\begin{equation*}
    u (t,\,0,\vec y) =2^{1-m}\pi^{-m/2}D^{1-m/2}\int\limits_0^t\int\limits_{\mathbb{R}^{m-1}} \exp\left[\frac{|\vec y-\vec y^{\:\prime} |^2}{4D(t-s)}\right]\partial_x u(s,\,0,\,\vec y^{\:\prime})
    \frac{\mathrm{d} \vec y^{\:\prime} \, \mathrm{d}s}{\left(t-s\right)^{m/2}}.
\end{equation*}

For the multidimensional problems the TBC for wave equations in the area $x>0$ at dimensions $m=2,\:3$ we have
\begin{equation*}
\partial_x u (t,\,0,\,y) = \frac{2c}{\sqrt{\pi}}\int\limits_0^t \int\limits_{|y-y^\prime|\le c(t-s)} 
\frac{\left[c^{-2}\partial_s^2- \partial_{y^\prime}^2\right] u(s,\,0,\,y^\prime)}{\sqrt{c^2 (t-s)^2 - |y-y^\prime|^2}} \, \mathrm{d}y^\prime \mathrm{d}s.
\end{equation*}
and
\begin{equation*}
\partial_x u (t,\,0,\,y,\,z) = \frac{2c}{\sqrt{\pi}}\int\limits_0^t \int\limits_{|y-y^\prime|^2+|z-z^\prime|^2\le c^2(t-s)^2} 
\left[c^{-2}\partial_s^2- \partial_{y^\prime}^2\right] u(s,\,0,\,y^\prime,\,z^\prime) \, \mathrm{d}y^\prime \mathrm{d}z^\prime \mathrm{d}s.
\end{equation*}

\textbf{Note 5}. {Sometimes the integration area in Eq.~(\ref{kernconvol2}) may be reduced. For instance,} in the case of hyperbolic equations in the half-space $\mathbb{R}^m_+$ the kernel's support is included into the inverse light cone $\left|\vec y-\vec z\right|\le c(t-s)$, where the constant $c$ is the speed of light, $m$ is odd, the kernel's support in these formulae may be reduced (by the Stokes formula) to the boundary of the inverse light cone $\left|\vec y-\vec z\right| = c(t-s)$, i.e. we obtain a lacuna (see for comparison \cite{atiyah1970lacunas}). The explicit formulae for TBC for other evolutionary partial derivative equations and systems of mathematical physics were considered in \cite{gordin1979diss}, see also  \cite{gordin2010mathematics}.

\textbf{Note 6}. Sometimes there is a finite-difference analogue of the inverse light cone for hyperbolic differential equations and Systs. (\ref{gener}). The DTBCs and {ADTBCs} were constructed in \cite{gordin1979diss} for explicit and “almost explicit” finite-difference schemes for the 2D wave equation, and the linearised shallow water system (barotropic model) with the Coriolis parameter. The DTBCs were constructed in the form of discrete convolutions with respect to time $t$ and 1D tangential (to the boundary) variable $y$. It is a scalar operator for the wave equation, and a matrix ($3\times 3$) operator for shallow water systems. The kernels of the convolutions were determined numerically. The inverse light cone was observed in both examples: the supports of the kernels are included in the cones, if the explicit schemes are stable (the Courant -- Friedrichs -- Lewy criterion is fulfilled).

Artificial boundary conditions for a {finite} area $V$ were studied using Calderon's integral operators in \cite{ryaben1990faithful, ryabenkii1993artificial, ryaben1995artificial, ryaben2006theoretical}. {DTBCs for Leontovich equation (which may be reduced to the Schr\"odinger equation) were constructed in \cite{fesch_2013} for cuboid $V$ by using the Fresnel transform.}

\section{Finite-Difference Implicit Scheme}
\label{se.schem}
Let us consider Eq.~(\ref{1.1}) with constant coefficients $\rho,\:R,\:E$ and with zeroth forcing $f$:
\begin{equation}
\label{2.1}     
\rho \, \frac{\partial^2 u}{\partial t^2} - R^2\rho \, \frac{\partial^4 u}{\partial t^2 \, \partial x^2}+ER^2 \, \frac{\partial^4 u}{\partial x^4}=0,
\end{equation}
and the implicit finite-difference scheme on the five-point stencil of the Crank -- Nicolson type  that approximates Eq.~(\ref{2.1}) on a uniform grid with the steps $\tau$ with respect to time $t$ and $h$ with respect to spatial variable $x$:

\begin{equation}
\label{2.2}     
    \sigma \left(u_{m+2}^{n+1} + u_{m-2}^{n+1} + u_{m+2}^{n-1} + u_{m-2}^{n-1}\right) + \beta \left(u_{m+1}^{n+1} + u_{m-1}^{n+1} + u_{m+1}^{n-1} + u_{m-1}^{n-1} \right) + \alpha \left(u_m^{n+1} + u_m^{n-1}\right) + \gamma \left(u_{m-1}^{n} + u_{m+1}^{n}\right) + \delta u_m^n = 0,
\end{equation}
where the upper index $n$ corresponds to the number of the time step, and the lower index $m$ corresponds to the spatial variable, $m=2,\,3,\,\ldots, N-2$, {where $N+1$ is the number of grid points in the segment $[-L/2,\,L/2]$}. Thus, we define $N-3$ linear algebraic equations with $N+1$ unknown values $\left\{u^{n+1}_m\right\}_{m=0}^N$. {The coefficients of the scheme are deduced using dimensionless parameters $\nu = E R^2 \rho^{-1} \cdot \tau^2 h^{-4}$, $\mu = R^2 \cdot h^{-2}$ and formulae $\alpha = 1 + 3\nu + 2\mu$, $\beta = -2\nu - \mu$, $\gamma = 2\mu$, $\delta = -2 - 4\mu$, $\sigma = \nu / 2$.}

To close the linear algebraic system we must add two linear algebraic equations that describe two boundary conditions on the left edge of the segment and, similarly, the same number of equations for the right edge. As a result, we obtain a {closed} linear algebraic system for $\left\{u^{n+1}_m\right\}_{m=0}^N$.

To begin the computational process and calculate the values $\left\{u^2_m\right\}_{m=0}^N$, we need two initial functions $\left\{u^0_m\right\}_{m=0}^N$ and $\left\{u^1_m\right\}_{m=0}^N$.

\section{{Discrete} Transparent Boundary Conditions, Rational Approximations, {Stability}}
\label{se.TBC}

\subsection{{{Derivation} of DTBCs for finite-difference equations}}
\label{subse.TBCsDef}
Consider an ordinary finite-difference equation of degree $n$ with constant coefficients
\begin{equation}
    \label{eq:OFDE}
    a_n u(m+n)+\ldots + a_0 u(m)=g(m),\quad m\in\mathbb{Z},
\end{equation}
where $g(m)$ is a given right hand side, which decreases as $m \to \pm \infty$. Further, consider the fundamental set of solutions as $g(m) \equiv 0$. If all roots of the characteristic equation
\begin{equation}
    \label{eq:OFDE-char}
    a_n \lambda^n + \ldots + a_0 = 0
\end{equation}
are different, then the fundamental set of solutions of Eq.~\eqref{eq:OFDE} can be expressed as $Y_k(m) = \mathrm{const}_k \cdot \lambda_k^m$, $k = 1,\ldots, n${, where $\lambda_k$ are the solutions of Eq.~\eqref{eq:OFDE-char}} . If some roots are multiple, we also have solutions $m^d \lambda_k^m$, where the degree $d$ is less than the multiplicity of the corresponding root $\lambda_k$. We also assume that {the boundary is not characteristic, i.e.} there are no roots of Eq.~\eqref{eq:OFDE-char} such that $|\lambda_k| = 1$.

Solutions of Eq.~\eqref{eq:OFDE} that are bounded as $m \to \infty$ could be expressed as (see, e.g. \cite{gelfond1971calculus})
\begin{equation}
    \label{eq:OFDE-sol}
    u(m) = \sum_{{j} = -\infty}^\infty G(m, j) g(j),
\end{equation}
where the kernel $G$ (Green function) is constructed using the fundamental set of solutions {on the real line.}

{Green's function (built by the fundamental set of solutions) is used to justify the algorithm of construction of DTBCs, see \cite{gordin1987mathematicalb, gordin1979diss, gordin2000mathematical}. For Schr\"odinger, wave, diffusion equations, and Eq.~\eqref{1.1} this fundamental set of solutions can be divided into two parts --- half of the solutions decreases as $x \to \infty$, the other half decreases as $x \to -\infty$. The similar statement is true for finite-difference equations, i.e. as $m \to \pm \infty$.}

Let $K$ be the number of roots of characteristic equation Eq.~\eqref{eq:OFDE-char} satisfying the inequality $|\lambda| < 1$, i.e.
\begin{equation*}
    |\lambda_1|\le\ldots \le |\lambda_K|<1< |\lambda_{K+1}|\le \ldots \le |\lambda_n|.
\end{equation*}

If we consider a partial differential equation, then the Fourier (and / or Laplace) transform should be applied to all variables that are tangential to the boundary (including time $t$), see \cite{gordin1978mixed}. For finite-difference equations, the analogous discrete transforms are used. The functions $Y_k$ are also dependent on dual variables. {We assume that the number $K$ does not depend on the dual variables, i.e. the border is not characteristic.} In Eq.~\eqref{2.2}, we only have one variable (time) that is tangential to the boundary.


If we consider the problem on the positive {discrete} half-line $m \ge 0$, then the summation in Eq.~\eqref{eq:OFDE-sol} is done only for non-negative indexes $j$. The kernel $G(m, j)$ remains the same (i.e. the solution of mixed problem for all decreasing functions $g(m)$ as $m \to \infty$), if at $m = 0$ there are such {$K$} boundary conditions that the functions of fundamental set of solutions are not changed. In other words, with such boundary conditions, the solution of the problem on the positive {discrete} half-line will coincide with the solution on the real line for $j \ge 0$. Therefore, such {$K$} boundary conditions imitate the bounded solution at $-\infty$, meaning that there is no reflection from the left boundary at $m=0$.

To obtain this transparency property {on the left border of Eq.~\eqref{eq:OFDE}}, it is sufficient {(see proof by author in \cite{gordin1978mixed} or \cite{gordin1979diss, gordin1987mathematical, gordin2000mathematical})} to construct homogeneous boundary conditions at point $m=0$, such that
\begin{enumerate}
    \label{TBCs-req}
    \item functions $Y_k(m)$ satisfy the conditions for all $k > K$;
    \item linear combination $\sum_{k=1}^K \mathrm{const}_k \, Y_k(m)$ satisfies the boundary conditions if and only if $\mathrm{const}_k = 0,\; \forall k = 1,\ldots,K$.
\end{enumerate}

{\textbf{Note 7.} To obtain the transparency property on the right border, the inequality sign in the requirement 1. should be changed to $k \leq K$, and the index $k$ in the requirement 2. (as well as the summation) is now from $K+1$ to $n$.}

The requirement 1. may be interpreted as an orthogonality condition in the space of linear boundary operators. The codimension of the subspace is equal to $n-K$.

\textbf{Note 8.} The requirement 2. is the discrete analogue of famous Shapiro -- Lopatinsky condition for differential equations. It can be formulated in the following form. Let us consider the matrix $\|I_{ij}\|_{i, j = 1}^K$, where $I_{ij}$ is the value of $i$-th boundary operator on $Y_j$. Its determinant (Lopatinsky determinant of the boundary problem) must be nonzero. The similar Lopatinsky determinant may be constructed for various differential problems (instead of the finite-difference).

{Later in Subsect.~\ref{subse.rational} we apply rational approximations of DTBCs to obtain ADTBC. Therefore, the corresponding Lopatinsky determinant will almost always be nonzero. To check the stability of mixed problem with applied ADTBCs, we propose other approach described in Subsect.~\ref{subsect.stability}.}

The requirements 1. and 2. do not define unique boundary conditions. The suitable gauge could be chosen in the space of boundary conditions -- in practice, specific boundary operators. We explain our choice of the gauge in Subsect.~\ref{subse.rational}.

If the problem is considered on a segment, then the boundary conditions that imitate the Cauchy problem are used on each edge. The total number of boundary conditions at both edges is equal to $n$. In this work, we have $n = 4$ and $K = 2$ (see Subsect.~\ref{subse.TBC}).

\subsection{General Plan of Approach}
\label{subse.algo}
The algorithm for constructing {ADTBCs}  at $x=\pm L/2$ for Eq. (\ref{2.2}) is as follows:

\begin{enumerate}[{\it Step 1.}]
\item
Apply the $\mathcal{Z}$-transformation (discrete analogue of the Laplace integral transformation) {with respect to time} to Eq.~(\ref{2.2}) and obtain a linear ordinary finite-difference equation with respect to the spatial variable $m$; the coefficients of the equation depend on the parameter $z\in \mathbb{C}$.     
 
\item 
Construct, for the corresponding homogeneous finite-difference 4-th order equation, a fundamental set of solutions $\left\{Y_j(m)\right\}_{j=1}^4$, such that  solutions $Y_1,\:Y_2$ decrease as $m\to +\infty$, and  solutions $Y_3,\:Y_4$ decrease as $m\to -\infty$.
 
\item 
Decompose the obtained growing (as $m\to +\infty$) solutions (functions from the dual variable $z$ with respect to time $t$) into a series in a neighbourhood of $z=\infty$.
 
 \item Construct vectorial rational functions (the construction {generalises} the Hermite -- Pad\'e approximation at the point $z=\infty$, see e.g. \cite{baker1996pade, gordin1982application, gordin1987mathematicalb, gordin2000mathematical, gordin2010mathematics}), which are asymptotically orthogonal to two growing solutions of this $\mathcal{Z}$-transformed equation. The corresponding polynomials are symbols of the transparent boundary operators.
\item
Apply the inverse $\mathcal{Z}$-transformation to the obtained coefficients of the {convolution operators in} {ADTBCs}.
\end{enumerate}

In order to calculate the inverse $\mathcal{Z}$-transformation, it is necessary to decompose the symbols of the corresponding operators into a Laurent series in the neighbourhood of the point $z=\infty$. For convenience, we introduce a change of variable: $\omega=1 / z$ and decompose the symbols into a Taylor series in the neighbourhood of the point $\omega=0$.

The growth and decrease of the Taylor coefficients of the meromorphic function are associated with the location of the function singularities on the complex plane. It is important to estimate the areas of the convergence of the obtained power series, which depend on the features of functions (solutions of the homogeneous finite-difference equation). {The singularities are either points of pole, or branch points.}

\subsection{Transparent Boundary Conditions for the Finite-Difference Equation}
\label{subse.TBC}
Let us apply the $\mathcal{Z}$-transformation with respect to time to Eq.~(\ref{2.2}). Then we obtain the linear ordinary non-homogeneous finite-difference equation
\begin{equation}
\label{4.1}     
    \sigma \left(z^2 + 1\right) \left[v(m+4) + v(m)\right] + \left(\beta \left( z^2 + 1 \right) + \gamma z \right) \left[v(m+3) + v(m + 1)\right] + \left( \alpha \left(z^2 + 1\right) + \delta z \right) v(m+2) = g(z, m),
\end{equation}
where $z\in\mathbb{C}$ is the variable that is dual to the discrete time $n$, $v(m) \equiv v(m,\,z)$ is the image of the solution $u^n_m$, and $g(z,m)$ is the right hand side that is obtained by the $\mathcal{Z}$-transformation from the initial functions $u^0_m$ and $u^1_m$. If we approximate non-homogeneous Eq.~~(\ref{1.1}), then the image of the $\mathcal{Z}$-transformation of the right hand side $f$ is included in $g (z,m)$.

The corresponding homogeneous equation after the change of variable $\omega=1/z$ has the form
\begin{equation}
\label{4.2}     
    \sigma \left(\omega^2 + 1\right) \left[v(m+4) + v(m)\right] + \left(\beta \left( \omega^2 + 1 \right) + \gamma \omega \right) \left[v(m+3) + v(m + 1)\right] + \left( \alpha \left(\omega^2 + 1\right) + \delta \omega \right) v(m+2) = 0.
\end{equation}

The order of the characteristic equation for ordinary finite-difference Eq.~(\ref{4.2}) (see e.g. \cite{gelfond1971calculus, gordin1978projectors, gordin1982boundary, gordin1987mathematicalb, gordin2000mathematical, gordin2010mathematics})
\begin{equation}
\label{4.3}     
    \sigma \left(1 + \omega^2\right) \left[\lambda^4 + 1\right] + \left(\beta \left(1 + \omega^2\right) + \gamma \omega \right) \left[\lambda^3 + \lambda \right] + \left(\alpha \left(1 + \omega^2\right) + \delta \omega \right) \lambda^2 = 0
\end{equation}
at $\omega\neq\pm i$ is equal to $4$ and is reciprocal. Let us divide Eq.~(\ref{4.3}) by $\lambda^2$ and rewrite it in the form
\begin{equation}
\label{4.4}
    \sigma \left(1 + \omega^2\right) \left[\lambda + \lambda^{-1}\right]^2 + \left(\beta \left(1 + \omega^2\right) + \gamma \omega \right) \left[\lambda + \lambda^{-1} \right] + \delta \omega + \left(\alpha -2\sigma \right) \left(1 + \omega^2\right) = 0.
\end{equation}

\textbf{Note 9.} The order of  algebraic Eq.~(\ref{4.3}) at $\omega = \pm i$ is equal to $3$:
\begin{equation*}
    \gamma\omega\left(\lambda^3+\lambda\right)+\delta\omega\lambda^2=0.
\end{equation*}

Therefore, $\lambda_1=0$, and as $\omega\to\pm i$ we obtain $\lambda_3\to\infty$. {Here the numbering of functions $\lambda_i$ is the same as in non-degenerate case Eq.~(\ref{4.3})}. We determine the remaining two roots of the cubic equation from the quadratic reciprocal equation:
\begin{equation*}
    \lambda^2 +\frac{\delta}{\gamma}\lambda + 1=0,
\end{equation*}
where $\delta/\gamma=-2-1/\mu$, and therefore 
\begin{equation*}
    \lambda_{2,4}=1+\frac{1}{2\mu}\mp\sqrt{\left(1+\frac{1}{2\mu}\right)^2-1}.
\end{equation*}

Thus, both roots are real and positive. According to Vieta's theorem, the following inequalities are fulfilled:
\begin{equation*}
    0 < \lambda_2 < 1 < \lambda_4.
\end{equation*}

For other values of $\omega$ we change variable $\eta=\lambda+\lambda^{-1}$ in Eq.~(\ref{4.4}) and obtain a quadratic equation for the auxiliary variable $\eta$:

\begin{equation}
\label{4.5}     
\sigma \left(1 + \omega^2\right) \eta^2 + 
\left(\beta \left(1 + \omega^2\right) + \gamma \omega \right) \eta + \delta \omega + \left(\alpha -2\sigma \right) \left(1 + \omega^2\right) = 0.
\end{equation}

The roots of Eq.~(\ref{4.5}) are
\begin{equation}
\label{4.6}   
\begin{aligned}
    \eta_1 (\omega) &= \frac{-\beta\left(1 + \omega^2 \right) - \gamma \omega - \sqrt[+]{\left(\beta \left(1 + \omega^2 \right) + \gamma \omega  \right)^2 - 4\sigma \left(1 + \omega^2 \right) \left[ \delta \omega + \left(\alpha - 2\sigma \right) \left(1 + \omega^2\right)\right]}}    {2\sigma \left(1 + \omega{^2} \right)}, \\
    \eta_2 (\omega) &= \frac{-\beta\left(1 + \omega^2 \right) - \gamma \omega + \sqrt[+]{\left(\beta \left(1 + \omega^2 \right) + \gamma \omega  \right)^2 - 4\sigma \left(1 + \omega^2 \right) \left[ \delta \omega + \left(\alpha - 2\sigma \right) \left(1 + \omega^2\right)\right]}}    {2\sigma \left(1 + \omega{^2} \right)},
\end{aligned}
\end{equation}
where $\sqrt[+]{y}$ is the complex root with a positive real part of $y\in \mathbb{C} \setminus \mathbb{R}_-$.

Before decomposing the functions $\lambda_j(\omega),\;j=\overline{1,\,4}$ into Taylor series in a vicinity of the origin $\omega=0$, we do it for the auxiliary functions $\eta_1(\omega),\:\eta_2(\omega)${, see \ref{app.eta12}:}
\begin{equation}
\label{eta1series}
    \eta_1(\omega) = \frac{1}{\nu} \sum_{k=0}^\infty (-1)^k \omega^{2k} \left[ \left(\mu+2\nu\right) \left(1+\omega^2\right) - 2\mu\omega - \sqrt{\mu^2-2\nu} \left(1-\omega\right) \left(\omega^2 - 2\frac{\mu^2}{\mu^2-2\nu} \omega + 1\right) \sum_{n=0}^\infty \mathrm{P}_n\left( \frac{\mu^2}{\mu^2-2\nu}\right) \omega^n \right],
\end{equation}
\begin{equation}
\label{eta2series}
\begin{aligned}
    \eta_2(\omega) = \frac{1}{\nu} \sum_{k=0}^\infty (-1)^k \omega^{2k} \left[ \left(\mu+2\nu\right) \left(1+\omega^2\right) - 2\mu\omega + \sqrt{\mu^2-2\nu} \left(1-\omega\right) \left(\omega^2 - 2\frac{\mu^2}{\mu^2-2\nu} \omega + 1\right) \sum_{n=0}^\infty \mathrm{P}_n\left( \frac{\mu^2}{\mu^2-2\nu}\right) \omega^n \right],
\end{aligned}    
\end{equation}
where $\mathrm{P}_n(\varepsilon)$ is a Legendre polynomial of degree $n$ at point $\varepsilon$. {The computation algorithm of Legendre polynomials is described in \ref{app.legendre}}.

Asymptotic of the functions as $\omega\to 0$ are described by the formulae
\begin{equation}
\label{4.10}   
\begin{aligned}
    \eta_1(\omega)=\vartheta_1 + r_1 (\omega),\\
\eta_2(\omega)=\vartheta_2 + r_2 (\omega),
\end{aligned}
\end{equation}
where $r_j(\omega)\to 0$, $j = 1, \, 2$, and

\begin{equation}
\label{4.doll}   
\begin{aligned}
    \vartheta_1 &= \frac{1}{2\sigma} \left[ \beta - \sqrt{\beta^2 - 4\sigma (\alpha - 2\sigma)} \right] = 2+\frac{\mu}{\nu} - \frac{1}{\nu} \sqrt{\mu^2-2\nu},\\
    \vartheta_2 &= \frac{1}{2\sigma} \left[ \beta + \sqrt{\beta^2 - 4\sigma (\alpha - 2\sigma)} \right] = 2+\frac{\mu}{\nu} + \frac{1}{\nu} \sqrt{\mu^2-2\nu}.
\end{aligned}
\end{equation}

If the inequality $\mu^2>2\nu$ is fulfilled, i.e. if
\begin{equation}
\label{4.vopr}   
\tau < R\sqrt{\frac{\rho}{2E}},
\end{equation}
then the radicand in Eq.~(\ref{4.doll}) is positive, and the values $\vartheta_1,\:\vartheta_2$ are real.

Let us resolve the relation $\eta=\lambda+\lambda^{-1}$ as a quadratic equation
\begin{equation}
    \label{4.9}
    \lambda^2-\eta\lambda+1=0.
\end{equation}

For both $\eta_1,\:\eta_2$ we obtain the following roots of characteristic Eq.~(\ref{4.3}):
\begin{equation*}
    \begin{aligned}
        \lambda_1 = \frac{\eta_1(\omega)}{2} - \sqrt{\frac{\eta_1^2(\omega)}{4} - 1},& \quad \lambda_2 = \frac{\eta_2(\omega)}{2} - \sqrt{\frac{\eta_2^2(\omega)}{4} - 1}, \\
        \lambda_3 = \frac{\eta_1(\omega)}{2} + \sqrt{\frac{\eta_1^2(\omega)}{4} - 1},& \quad \lambda_4 = \frac{\eta_2(\omega)}{2} + \sqrt{\frac{\eta_2^2(\omega)}{4} - 1}.
    \end{aligned}
\end{equation*}

\textbf{Note 10}. According to Vieta's theorem, either the absolute values of both roots of Eq.~(\ref{4.9}) are equal to $1$, or one absolute value is smaller than $1$ and the other one is greater. The first version takes place{, if} $\eta \in \left[-2, \, 2\right] \subset \mathbb{R}$. As $\omega \to 0 \in \mathbb{C}$ {this can never be obtained. Indeed}, if inequality (\ref{4.vopr}) is fulfilled, the square roots in Eq.~(\ref{4.doll}) are real and $2<\vartheta_1<\vartheta_2$.
If the equality $\mu^2=2\nu$ is fulfilled, then $2<\vartheta_1=\vartheta_2$.
If the inequality, which {is} inverse to (\ref{4.vopr}) is fulfilled, then the values $\vartheta_1,\:\vartheta_2$ are not real.
Thus, at small $\omega$ we obtain $\eta_{1,2}\notin [-2,\,2]$, i.e. the necessary condition of the mixed initial-boundary problem correctness for the finite-difference Eq.~(\ref{2.2}) is fulfilled.

The Taylor series of functions $\lambda_i(\omega)$ at point $\omega = 0$ have the forms {(see \ref{app.lambda})}

\begin{equation}
\label{4.11}
    \lambda_{1,3}(\omega) = \frac{\eta_1(\omega)}{2} \mp \sqrt{\frac{\vartheta_1^2}{4} - 1} \; \cdot \sum_{n=0}^\infty \frac{(-1)^n \, (2n)! \, r_1^n(\omega)}{(1-2n) \, n! \, 4^n \, (\theta_1+2)^n} \cdot \sum_{n=0}^\infty \frac{(-1)^n \, (2n)! \, r_1^n(\omega)}{(1-2n) \, n! \, 4^n \, (\theta_1-2)^n},
\end{equation}
\begin{equation}
\label{4.12}
    \lambda_{2,4}(\omega) = \frac{\eta_2(\omega)}{2} \mp \sqrt{\frac{\vartheta_2^2}{4} - 1} \; \cdot \sum_{n=0}^\infty \frac{(-1)^n \, (2n)! \, r_2^n(\omega)}{(1-2n) \, n! \, 4^n \, (\theta_2+2)^n} \cdot \sum_{n=0}^\infty \frac{(-1)^n \, (2n)! \, r_2^n(\omega)}{(1-2n) \, n! \, 4^n \, (\theta_2-2)^n},
\end{equation}
where $\eta_1(\omega)$, $\eta_2(\omega)$ and $r_{1,2}(\omega)$ are taken from Eqs.~\eqref{eta1series}, \eqref{eta2series} and \eqref{4.10},  {respectively}.

The following inequalities are fulfilled as $\omega \to 0$
\begin{equation*}
    |\lambda_1|, \, |\lambda_2| < 1 < |\lambda_3|, \, |\lambda_4|.
\end{equation*}

Therefore, as $m\to +\infty$ it is possible to derive decreasing $\lambda_1^m$, $\lambda_2^m$ and increasing $\lambda_3^m$, $\lambda_4^m$ solutions of finite-difference Eq.~(\ref{4.2}), which form the fundamental set of solutions.

\subsection{Transparent Boundary Conditions and Rational Approximations}
\label{subse.rational}

As with differential Eq.~(\ref{2.1}), for the correctness of the mixed initial-boundary value problem for finite-difference Eq.~(\ref{2.2}), two boundary conditions at each edge of the rod are required. Thus, the values of the solution $u$ at boundary and {several preboundary grid points} should be calculated at every time step. We start by constructing the $\mathcal{Z}$-image of the boundary conditions for the left edge in the form:
\begin{equation*}
\label{5.1}
\begin{aligned}
    P_1(\omega)\,v(0) + Q_1(\omega)\,v(1) &+ R_1(\omega)\,v(2) + S_1(\omega)\,v(3) = 0, \\
    P_2(\omega)\,v(0) + Q_2(\omega)\,v(1) &+ R_2(\omega)\,v(2) + S_2(\omega)\,v(3) = 0, 
\end{aligned}
\end{equation*}
which (after the inverse $\mathcal{Z}$-transformation) corresponds to the relations
\begin{equation}
\label{5.star}
    \sum_{j=0}^\infty p_{kj}\,u_0^{n-j} + \sum_{j=0}^\infty q_{kj}\,u_1^{n-j} + \sum_{j=0}^\infty r_{kj}\,u_{2}^{n-j} + \sum_{j=0}^\infty s_{kj}\,u_{3}^{n-j} = 0, \quad k=1,\,2,
\end{equation}
where values $p_{kj}$, $q_{kj}$, $q_{kj}$ and $r_{kj}$ are the coefficients in Taylor series before the term $\omega^j$ of the functions $P_k(\omega)$, $Q_k(\omega)$, $R_k(\omega)$, $S_k(\omega)$,  {respectively} ($k = 1, \, 2$) in a vicinity of the point $\omega=0\in\mathbb{C}$. {The requirements 1. and 2. from Subsect.~\ref{subse.TBCsDef} must be fulfilled for all $\omega$ in the vicinity of zero.}

Two (with numbers $k=1,\,2$) linearly independent boundary conditions will provide the transparency property, {if and only if} for the increasing Cauchy problem solutions $\nu (m)=\lambda_3^m$ and $\nu (m)=\lambda_4^m$ the symbols of the boundary conditions $\langle P_1,\,Q_1,\,R_1,\,S_1\rangle$ and  $\langle P_2,\,Q_2,\,R_2,\,S_2\rangle$ fulfill the following equations:
\begin{equation}
\label{and}
\begin{aligned}
    P_k + Q_k \, \lambda_3 + & R_k \, \lambda_3^2 + S_k \, \lambda_3^3 = 0, \\
    P_k + Q_k \, \lambda_4 + & R_k \, \lambda_4^2 + S_k \, \lambda_4^3 = 0.
\end{aligned}
\end{equation}

{\textbf{Note 11}. Here we compose spatial stencils for the boundary conditions from $4$ points with respect to $m$: $-L/2$, $-L/2+h$, $-L/2+2h$, $-L/2+3h$. However, more space steps could be included:}

{
\begin{equation*}
    \begin{aligned}
    \sum_{i = 0}^I P_{k, i} \, \lambda_3^i &= 0, \\
    \sum_{i = 0}^I P_{k, i} \, \lambda_4^i &= 0.
\end{aligned}
\end{equation*}}

{
Here, $I$ is a number of space steps used in the boundary conditions, $P_{k, i} \equiv P_{k, i}(\omega)$ are polynomials, index $k$ corresponds to number of boundary condition ($k = 1$ for the first boundary condition and $k=2$ for the second one), and index $i$ corresponds to the space step ($i = 0$ for boundary layer, $i = 1$ for preboundary layer, etc).}

For any given values $\omega$ the subspace of solutions of  Syst.~(\ref{and}) is two-dimensional and two boundary conditions on every edge can provide a boundary problem correctness for finite-difference Eq.~(\ref{4.1}). However, the subspace of the solutions of  Syst.~(\ref{and}) in the space of analytic vector-functions of $\omega$ is infinity-dimensional. 

The meromorphic functions $\lambda_3(\omega),\: \lambda_4(\omega)$ are not rational, and usually solutions of Syst.~(\ref{and}) in the subspace of polynomials do not exist. Therefore, DTBCs are non-local with respect to time. They include values of the solution of the boundary problem for Eq.~(\ref{2.2}) in the infinite number of previous time moments. For such a realisation of the {ADTBCs} the number of arithmetic operations, as well as the necessary computer memory, is proportional to the number of the temporal steps $n$. {We say that the boundary conditions are local with respect to time, if there exists a fixed natural number $n_t$ (number of time steps), such that at any time step $T = n \tau, \, n \in \mathbb{N}$ the boundary value is calculated only by using several preboundary points at the times $n\tau$, $(n-1)\tau$, $(n-2)\tau$, $\ldots$, $(n-n_t+1)\tau$.}

That is why we relax the requirements for the symbols of the operators of DTBCs, and exchange unknown analytic functions in Syst.~(\ref{and}) by polynomials and exact equalities by asymptotic (as $\omega\to 0$ and $k=1,\,2$):
\begin{equation}
\label{5.and}
\begin{cases}
    P_k(\omega) + Q_k(\omega)\,\lambda_3(\omega) + R_k(\omega)\,\lambda_3^2(\omega) + S_k(\omega)\,\lambda_3^3(\omega) = O\left(\omega^{K_k}\right), \\
    P_k(\omega) + Q_k(\omega)\,\lambda_4(\omega) + R_k(\omega)\,\lambda_4^2(\omega) + S_k(\omega)\,\lambda_4^3(\omega) = O\left(\omega^{K_k}\right).
\end{cases}
\end{equation}

The ``physical'' interpretation of the exchange: we neglect the impact of solution's values {for the distant}  past, assuming the resulting opacity is small. {However the requirement 2. from Subsect.~\ref{subse.TBCsDef} may be violated as the result of the exchange of Eq.~\eqref{and} to Eq.~\eqref{5.and}. Then, the corresponding mixed problem will be incorrect.}

We fix the degrees of the polynomials $P_k$, $Q_k$, $R_k$, and $S_k$, i.e. stencils for the {ADTBCs}. The number of the degrees of freedom in these four polynomials is equal to $M_k=\mathrm{deg}\,P_k+\mathrm{deg}\,Q_k+\mathrm{deg}\,R_k+\mathrm{deg}\,S_k+4$, and the number of  linear algebraic equations, which are obtained from Syst.~(\ref{5.and}) is equal to $2K_k$. Together with two normalisation conditions, which will be considered below, we obtain $2K_k+2$ linear algebraic equations. Thus, if $M_k=2K_k + 2$ and the system of linear algebraic equations is non-degenerate, we determine a unique set of the polynomials $P_k$, $Q_k$, $R_k$, and $S_k$ with given degrees.

If we choose normalisation condition at $k=1$:
\begin{equation*}
    P_1(0) = p_{1,0}=1, \quad Q_1(0) =q_{1,0}= 0,
\end{equation*}
then we are able to compute the value $u^n_0$ for every temporal step $n$ using the values in the inner points $u^j_{2h},\: u^j_{3h}$ at $j\le n$ and in the border and preborder points at previous time moments: $u^j_{0},\: u^j_{1}$ at $j< n$. Thus, we obtained the first boundary condition in the form:
\begin{equation}
\label{2zv}
 u_0^n + \sum_{j=1}^{\mathrm{deg}\,P_1} p_{1j}\,u_0^{n-j} + \sum_{j=1}^{\mathrm{deg}\,Q_1} q_{1j}\,u_1^{n-j} + \sum_{j=0}^{\mathrm{deg}\,R_1} r_{1j}\,u_{2}^{n-j} + \sum_{j=0}^{\mathrm{deg}\,S_1} s_{1j}\,u_{3}^{n-j} = 0.
\end{equation}

If we choose normalisation condition at $k=2$:
\begin{equation*}
    P_2(0) = p_{2,0}=0, \quad Q_2(0) =q_{2,0}= 1,
\end{equation*}
then we are able to compute the value $u^n_1$ for every temporal step $n$ using the values in the inner points $u^j_{2},\: u^j_{3}$ at $j\le n$ and in the border and preborder points at previous time moments: $u^j_{0},\: u^j_{1}$ at $j< n$. So, we obtain the second {ADTBCs}:
\begin{equation}
\label{3zv}
 u_1^n+   \sum_{j=1}^{\mathrm{deg}\,P_2} p_{2j}\,u_0^{n-j} + \sum_{j=1}^{\mathrm{deg}\,Q_2} q_{2j}\,u_1^{n-j} + \sum_{j=0}^{\mathrm{deg}\,R_2} r_{2j}\,u_{2}^{n-j} + \sum_{j=0}^{\mathrm{deg}\,S_2} s_{2j}\,u_{3}^{n-j} = 0.
\end{equation}

\textbf{Note 12.} Conditions in Syst.~(\ref{5.and}) are similar to the famous Hermite -- Pad\'e rational approximation of meromorphic functions. However, it is identical in the case of unique {ADTBCs} (see \cite{gordin1982application, gordin1987mathematicalb, gordin2000mathematical, gordin2010mathematics}). For two {ADTBCs} it is a more general vector rational approximation. Similar constructions, which can be considered as a generalisation of the Hermite -- Pad\'e rational approximation, were studied in \cite{nuttall1984asymptotics}.

The coefficients of Syst.~(\ref{5.and}) are not real, because the functions $\lambda_3(\omega)$ and $\lambda_4(\omega)$ are, in general, complex. As a result, the characteristic values are complex conjugated, and we can transform the systems to the real form:

\begin{equation}
\label{5.4}
\begin{cases}
    P_1(\omega)+Q_1(\omega)\,\frac{\lambda_3(\omega)+\lambda_4(\omega)}{2} + R_1(\omega)\,\frac{\lambda_3^2(\omega)+\lambda_4^2(\omega)}{2} + S_1(\omega)\,\frac{\lambda_3^3(\omega)+\lambda_4^3(\omega)}{2} = O \left(\omega^{K_1}\right), \\
    Q_1(\omega)\,\frac{\lambda_3(\omega)-\lambda_4(\omega)}{2\text{i}} + R_1(\omega)\,\frac{\lambda_3^2(\omega)-\lambda_4^2(\omega)}{2\text{i}} + S_1(\omega)\,\frac{\lambda_3^3(\omega)-\lambda_4^3(\omega)}{2\text{i}} = O \left(\omega^{K_1}\right), \\
    P_1(0) = 1, Q_1(0) = 0,
\end{cases}
\end{equation}
at $k=1$, and
\begin{equation}
\label{5.5}
\begin{cases}
    P_2(\omega)+Q_2(\omega)\,\frac{\lambda_3(\omega)+\lambda_4(\omega)}{2} + R_2(\omega)\,\frac{\lambda_3^2(\omega)+\lambda_4^2(\omega)}{2} + S_2(\omega)\,\frac{\lambda_3^3(\omega)+\lambda_4^3(\omega)}{2} = O \left(\omega^{K_2}\right), \\
    Q_2(\omega)\,\frac{\lambda_3(\omega)-\lambda_4(\omega)}{2\text{i}} + R_2(\omega)\,\frac{\lambda_3^2(\omega)-\lambda_4^2(\omega)}{2\text{i}} + S_2(\omega)\,\frac{\lambda_3^3(\omega)-\lambda_4^3(\omega)}{2\text{i}} = O \left(\omega^{K_2}\right), \\
    P_2(0) = 0, Q_2(0) = 1,
\end{cases}
\end{equation}
at $k=2$.

We solve Syst.~\eqref{5.4} and Syst.~\eqref{5.5} separately and find the coefficients of the {ADTBCs} on the left edge of the rod. Using the same approach, we find coefficients for the {ADTBCs} on the right edge. The results of numerical experiments with such {ADTBCs} are presented in Sect.~\ref{se.resul}.

\subsection{Stability of mixed problem with {ADTBCs}}
\label{subsect.stability}

The correctness (stability) of the Cauchy problem for a differential or finite difference equation does not guarantee that the mixed boundary value problem will also be correct. For stability, the boundary should be uncharacteristic and the boundary operators should satisfy the Shapiro -- Lopatinsky conditions, see \cite{agranovich1971boundary, sakamoto1970ii, sakamoto1970i, Kreiss1970, gordin1987mathematical, gordin2000mathematical, gordin2010mathematics, gordin1987mathematicalb, KreissHedwig2001, hormander2007analysis}. If the DTBCs were perfectly realised according to Eq.~\eqref{and}, the solution of the mixed problem would coincide with the solution of the Cauchy problem, i.e. the problem would be stable. However, we are forced to implement the DTBCs according to Syst.~\eqref{5.4} and Syst.~\eqref{5.5}, i.e. approximately. Thus, we need to check the stability of the mixed boundary value problem again. {Here, we propose a numerical method to do so.}

The Crank -- Nicolson approximation Eq.~\eqref{2.2} of the Cauchy problem for rod transverse vibrations Eq.~\eqref{2.1} is absolutely stable (see \ref{app_CN_stability}).

Usually, Chebyshev's norm in the space of grid functions is used to assess the stability of a finite-difference scheme:
\begin{equation}
\label{norm:chebyshev}
    \|u^n\|_{\mathbf{C}} = \max_{0 \leq j \leq N} |u^n_j|,
\end{equation}
where $u_j^n$ is a grid function, $n$ is the number of time step, $j$ is the number of space step. Also, the approximation of $\mathbf{L}^2$ norm is used (trapezoidal approximation):
\begin{equation}
\label{norm:L2}
    \|u^n\|_{\mathbf{L}^2} = \sqrt{h \left[ \frac{1}{2} \left((u_0^n)^2 + (u_N^n)^2\right) + \sum_{j = 1}^{N-1} (u_j^n)^2 \right]}.
\end{equation}

In the case of Eq.~\eqref{1.1}, the norms \eqref{norm:chebyshev} and \eqref{norm:L2} of Cauchy problem's solution (in an infinite domain) may increase over time. At the same time, the energy norm \eqref{norm:H} is preserved. Hence, the dynamic of a rod's energy is the key factor in the  estimation of stability (finite-difference approximation and boundary conditions). For the finite domain $x \in [-L/2, \, L/2]$ for the model with {ADTBCs} at any time moment $t>0$ this energy should be less than at time moment $t=0$. This means that a part of the energy was transferred from the segment outside, i.e. there was no reflection from the boundaries.

In our experiments we use an approximation of Hamiltonian $\mathcal{H}$, see Eq.~\eqref{eq:kinetic_energy} and \eqref{eq:potential_energy}, at the time moment $\tau(n + 1/2)$:
\begin{equation}
    \label{eq:hamiltonian_approx}
    \hat{\mathcal{H}}\left[u^{n+1/2}\right] = h \left[ \frac{1}{2} \left( \vartheta_0^{n+1/2} + \vartheta_N^{n+1/2} \right) + \sum_{j = 1}^{N-1} \vartheta_j^{n+1/2} \right],\quad \|u^{n+1/2}\|_\mathcal{H}=\sqrt{\hat{\mathcal{H}}\left[u^{n+1/2}\right],}
\end{equation}
where for $j = 1, \ldots, N-1$ we have
\begin{equation*}
    \vartheta_j^{n+1/2} =\rho \left(\frac{u_j^{n+1}-u^n_j}{\tau}\right)^2 + \rho R^2 \left(\frac{u_{j+1}^{n+1}-u^n_{j+1}-u_{j-1}^{n+1}+u^n_{j-1}}{2h\tau}\right)^2+ER^2\left(\frac{u_{j+1}^{n+1}-2u^{n+1}_{j}+ u^{n+1}_{j-1}+u_{j+1}^{n}-2u^{n}_{j}+ u^{n}_{j-1}}{2h^2}\right)^2,
\end{equation*}
and for the boundary values at $j = 0$ and $j = N$ the derivatives with respect to space $x$ and time $t$ should be approximated with forward (or backward) differences.

As stated in Subsect.~\ref{subse.rational}, we can choose any set of polynomial degrees in Syst.~\eqref{5.and}. However, the existence of the system's solution (with respect to unknown polynomial coefficients) does not guarantee the stability of the mixed initial-boundary problem's \eqref{2.2} solution. {ADTBCs} could result in a partial reflection of outgoing waves back to the computational domain, solution's ``explosion'' (i.e. instability), or Syst.~\eqref{5.and} could not be solved for a given set of polynomial degrees (i.e. {ADTBCs} do not exist).

If the set of polynomial coefficients does indeed provide the transparency of the border, then the corresponding solution should decrease over time. In fact, the waves spread outside of the computational domain, see Fig.~\ref{fig:2}. Therefore, we say that a finite-difference scheme with {ADTBCs} is stable, if at any time moment the energy norm of the solution obtained with these {ADTBCs} is less or equal to the energy norm at the initial time moment $t = 0$, i.e.
\begin{equation}
\label{eq:tbcs_stable}
\sqrt{\hat{\mathcal{H}}\left[u^{n+1/2}\right]} \leq \sqrt{\hat{\mathcal{H}}\left[u^{1/2}\right]} \quad \forall n = 0, \ldots, N_t-1,
\end{equation}
where we set the number of time steps to $N_t=10^5$. Along with energy criterion \eqref{eq:tbcs_stable} we introduce the classical stability criteria of $\mathbf{C}$-norm
\begin{equation}
    \label{crit:C}
    \|u^n\|_{\mathbf{C}} \leq \|u^0\|_{\mathbf{C}}, \quad \forall n = 0, \ldots, N_t,
\end{equation}
and $\mathbf{L}^2$-norm
\begin{equation}
    \label{crit:L2}
    \|u^n\|_{\mathbf{L}^2} \leq \|u^0\|_{\mathbf{L}^2}, \quad \forall n = 0, \ldots, N_t.
\end{equation}

\begin{figure}[ht]
    \centering
    \begin{subfigure}[b]{.33\linewidth}
        \includegraphics[width = \textwidth]{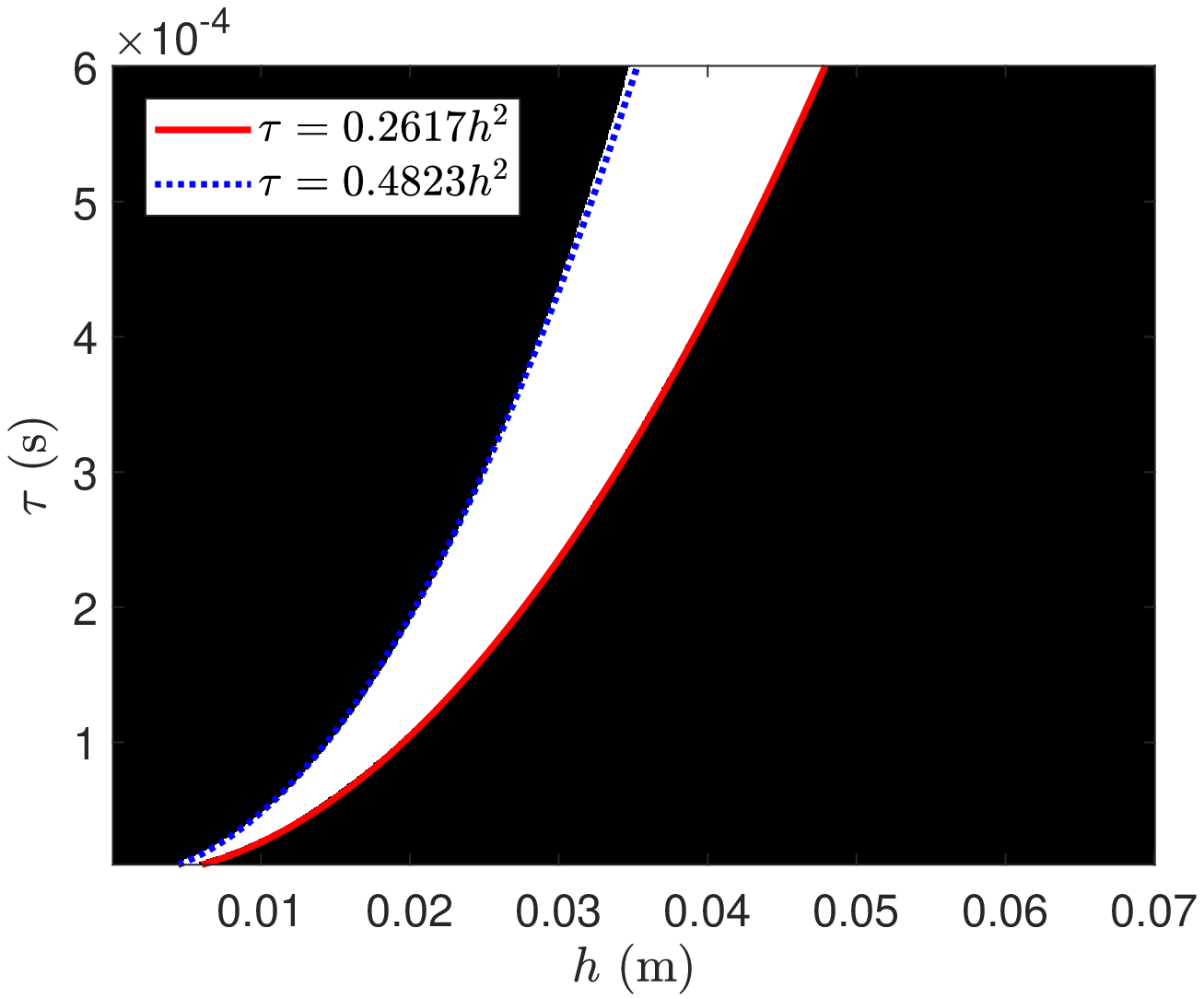}
        \caption{\label{fig:1a}}
    \end{subfigure}
    \begin{subfigure}[b]{.33\linewidth}
        \includegraphics[width = \textwidth]{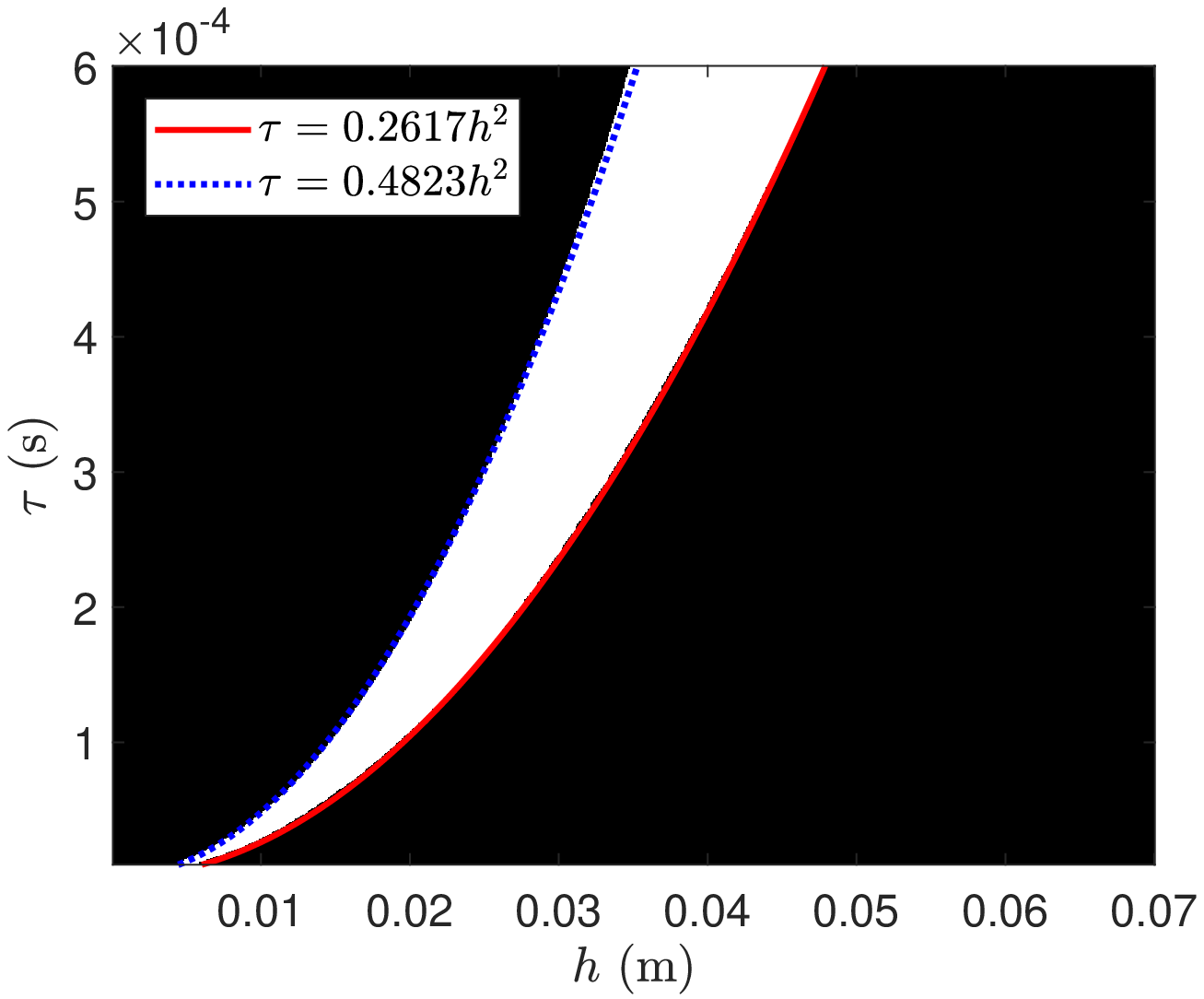}
        \caption{\label{fig:1b}}
    \end{subfigure}
    \begin{subfigure}[b]{.33\linewidth}
        \includegraphics[width = \textwidth]{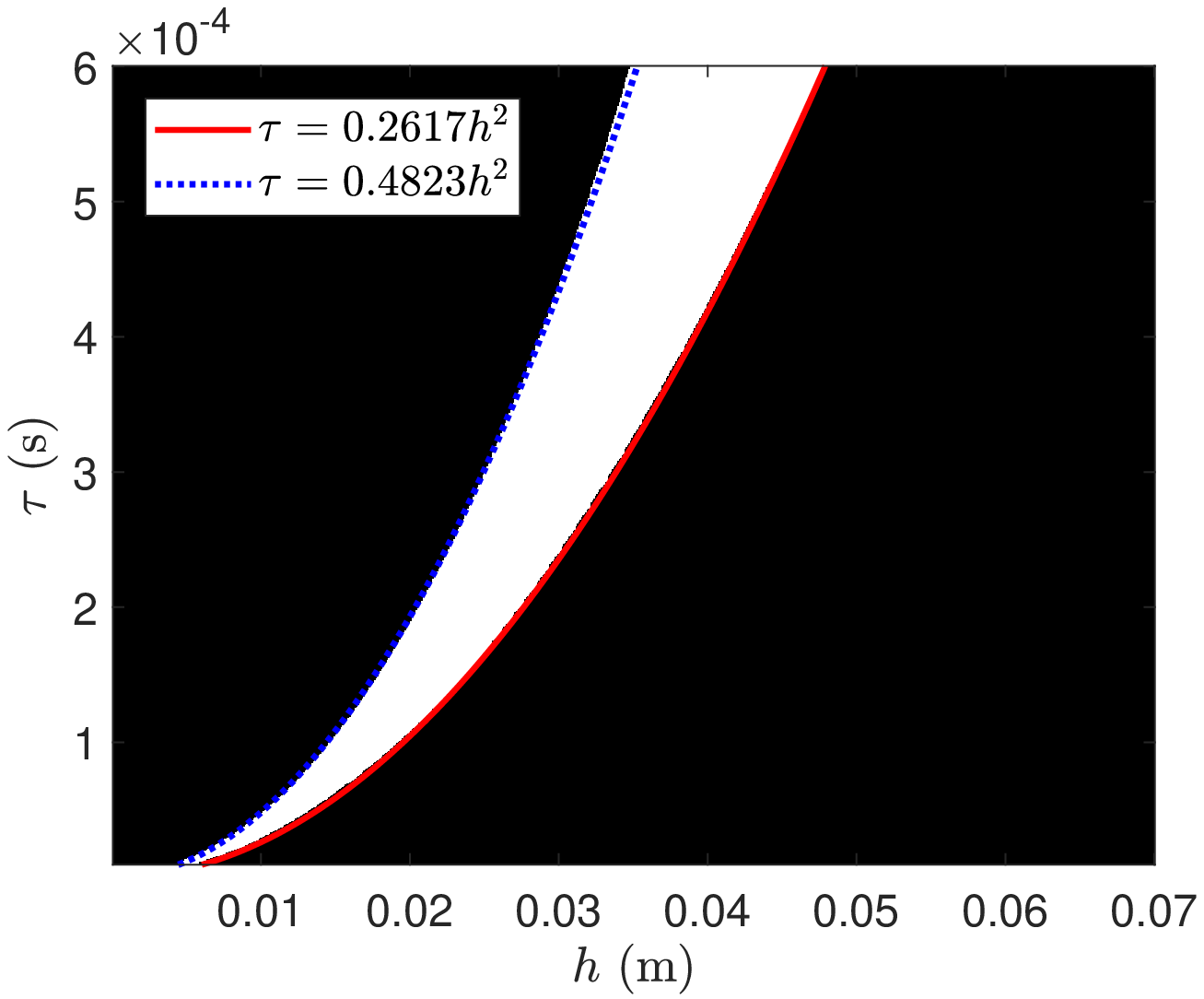}
        \caption{\label{fig:1c}}
    \end{subfigure}
    \begin{subfigure}[b]{.33\linewidth}
        \includegraphics[width = \textwidth]{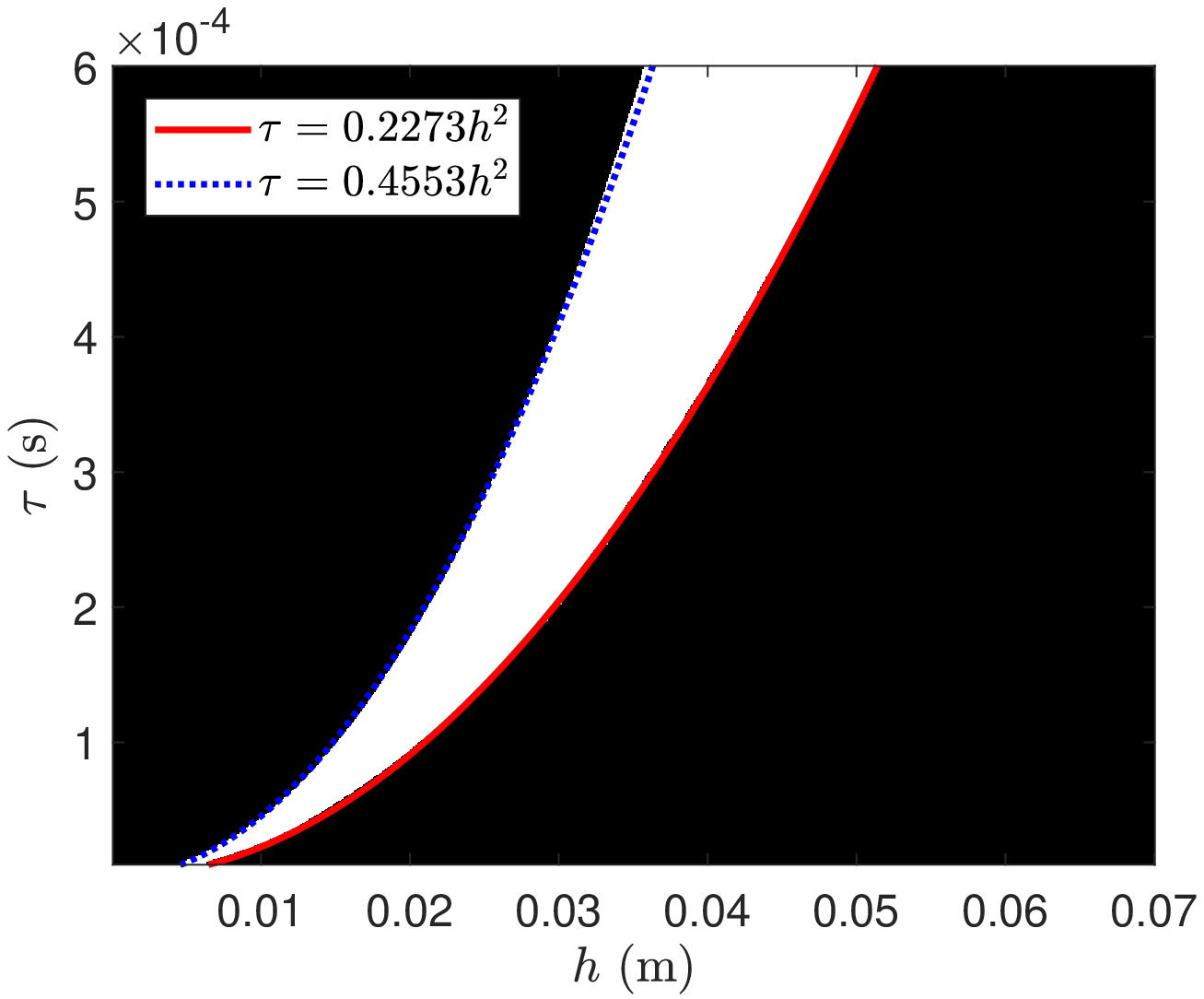}
        \caption{\label{fig:1d}}
    \end{subfigure}
    \begin{subfigure}[b]{.33\linewidth}
        \includegraphics[width = \textwidth]{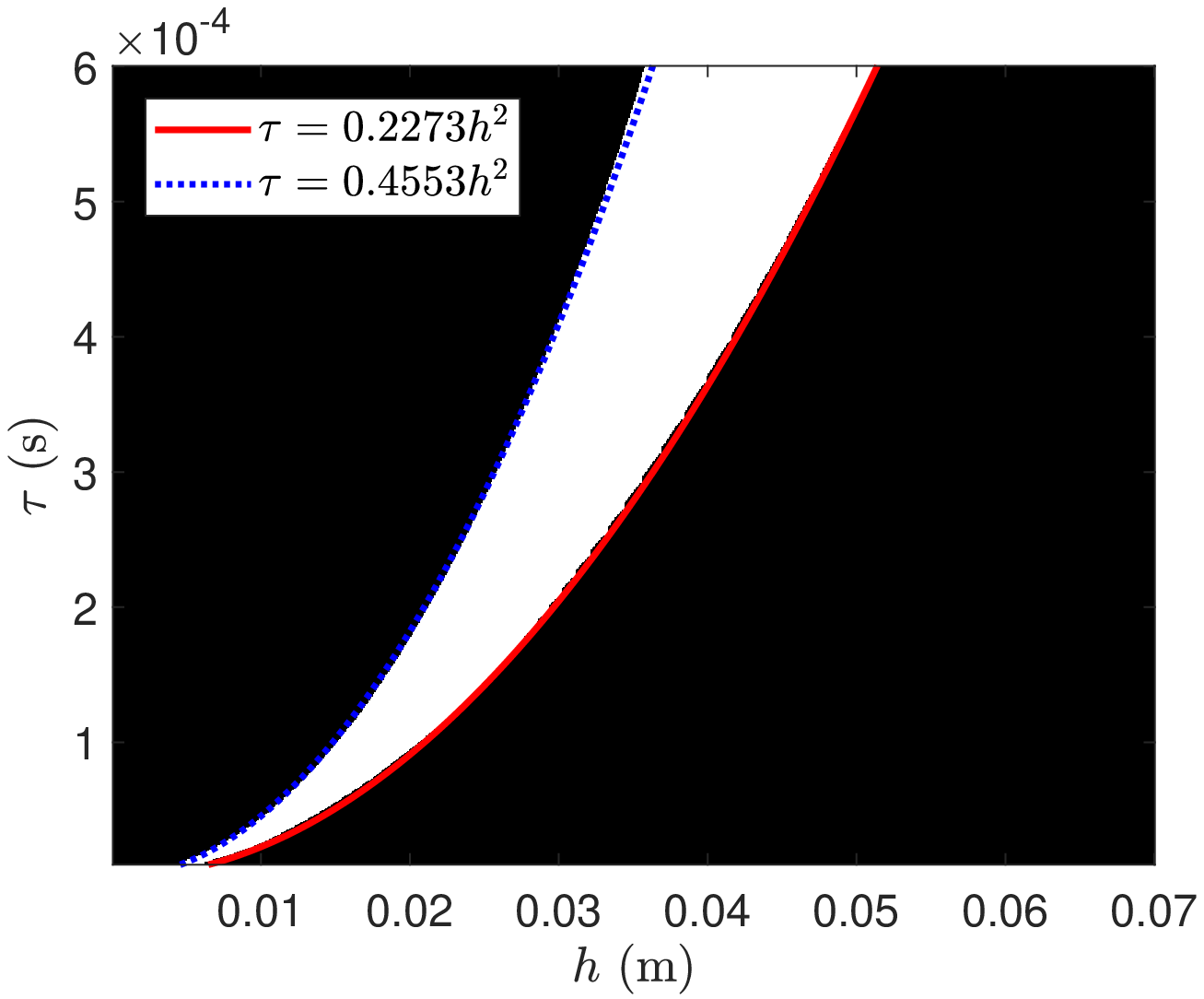}
        \caption{\label{fig:1e}}
    \end{subfigure}
    \begin{subfigure}[b]{.33\linewidth}
        \includegraphics[width = \textwidth]{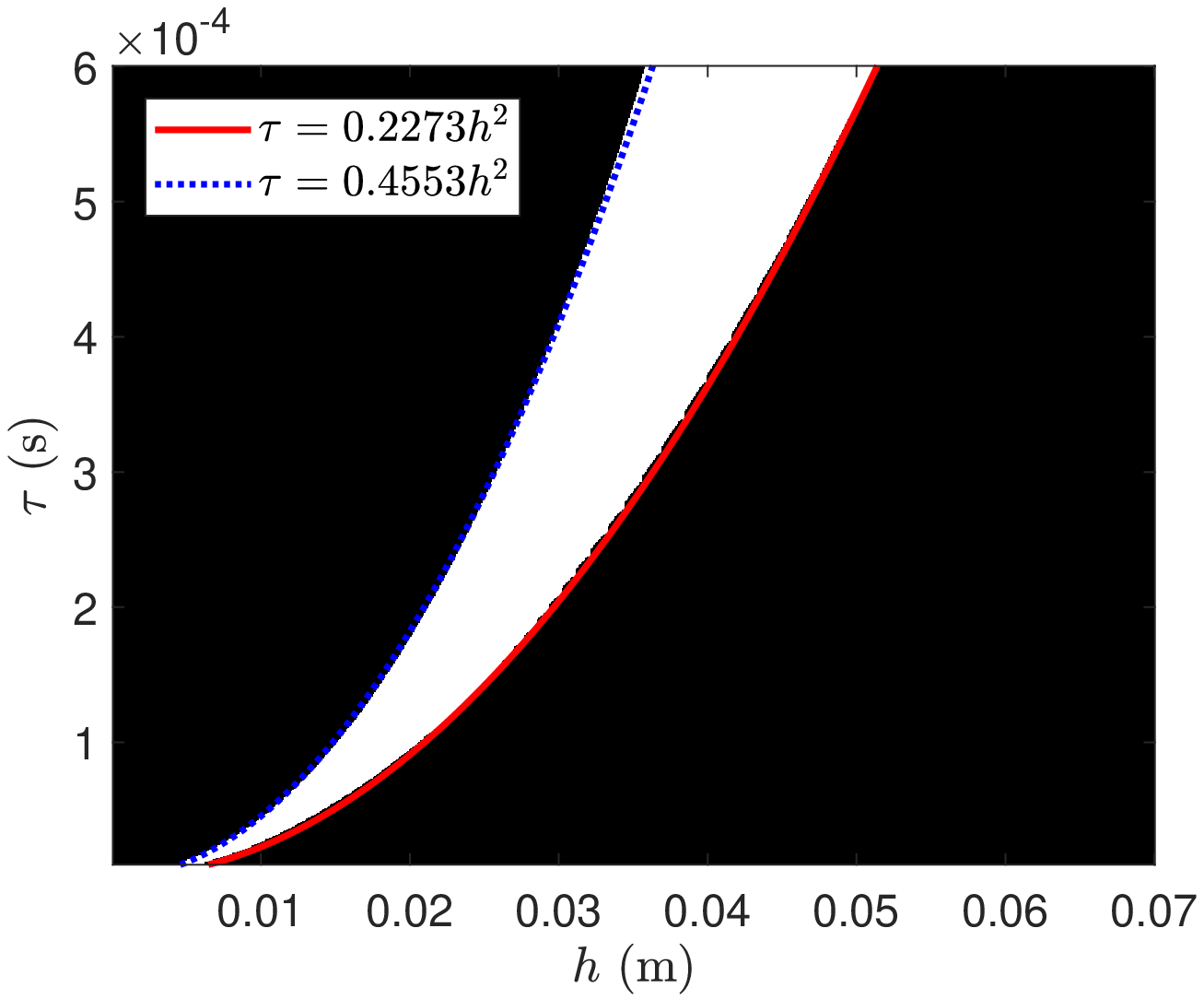}
        \caption{\label{fig:1f}}
    \end{subfigure}
    \caption{The domain of stability on the $(h, \tau)$ plane for two {ADTBCs}. Left -- energy norm criterion, centre -- $\mathbf{C}$-norm criterion, right -- $\mathbf{L}^2$-norm criterion. The rows correspond to the particular set of polynomial degrees: (a), (b), (c) -- $\langle 4, 4, 8, 8 \rangle$; (d), (e), (f) -- $\langle 5, 3, 9, 7 \rangle$. The white area -- stable, conditions \eqref{eq:tbcs_stable}, \eqref{crit:C} and \eqref{crit:L2} are fulfilled (energy, $\mathbf{C}$ and $\mathbf{L}^2$ norms, respectively). The black area -- stability conditions are not fulfilled. {Physical parameters of the rod $\rho$, $E$, $R$ and $L$ are the same as in Table~\ref{tab:model_params}.
    }}
    \label{fig:1}
\end{figure}

The domains of stability (white) are presented in Fig.~\ref{fig:1}. For each case we approximate the boundaries between stable and unstable regions.

Unlike the classical equations of mathematical physics (wave, heat, Schr\"odinger equations), the classical stability criteria \eqref{crit:C} and \eqref{crit:L2} for Eq.\eqref{2.1} are not necessarily fulfilled. 

In our experiments for simplicity we used $U_1 \equiv 0$. For finite-difference Eq.~\eqref{2.2} all these {criteria} are similar, see Fig.~\ref{fig:1}. {There are small differences between all norms that occur near the boundaries between white and black regions. However, the approximation parabolas remain unchanged throughout all criteria, see Fig.~\ref{fig:1a}, \ref{fig:1b}, \ref{fig:1c} (Fig.~\ref{fig:1d}, \ref{fig:1e}, \ref{fig:1f}).}

We denote the symbol of {ADTBCs} obtained with polynomial degrees $\mathbf{deg} P_k = d_{1,k}$, $\mathbf{deg} Q_k = d_{2,k}$, $\mathbf{deg} R_k = d_{3,k}$ and $\mathbf{deg} S_k = d_{4,k}$ with $k = 1,\,2$ as
\begin{equation*}
    \langle P_k, Q_k, R_k, S_k \rangle \equiv \langle d_{1,k}, d_{2,k}, d_{3,k}, d_{4,k} \rangle.
\end{equation*}

Eq.~\eqref{2.1} {on the segment [-L/2, \, L/2]} could be rewritten as
\begin{equation}
\label{eq:rod_rewritten}
    \frac{\partial^2 u}{\partial t^2} - D \frac{\partial^4 u}{\partial t^2\partial x^2} + C\frac{\partial^4 u}{\partial x^4}=0,
\end{equation}
where the physical dimensions of the constants $L$, $D=R^2$, and $C=ER^2\rho^{-1}$ are $\mathrm{m}$, $\mathrm{m}^2$, and $\mathrm{m}^4\,\mathrm{s}^{-2}$, respectively.

{Usually, the stability condition of a finite-difference scheme of the Eq.~\eqref{2.2} type has the form $\tau < A h^2$, where the constant $A$ has the physical dimension $\mathrm{m}^{-2}\,\mathrm{s}$ and depends on physical parameters and boundary conditions. Here, the time dimension is present only in the constant $C=ER^2\rho^{-1}$. Dimensionless parameter can be obtained only as a function of $\omega = L/R$. Therefore, from the dimension theory \cite{barenblatt1996scaling}, it follows that $A = B(\omega) C^{-1/2}$.}

{From our numerical experiments, surprisingly, the time step $\tau$ should be bounded from both sides, which does not contradict the dimension theory. Moreover,} the stability domain could be composed from several parabolic sectors:
\begin{equation*}
    A_1 h^2 < \tau < A_2 h^2,
\end{equation*}
where both constants could be expressed as
\begin{equation*}
    A_1=B_1 (\omega) \, C^{-1/2},\quad A_2=B_2(\omega) \, C^{-1/2},\quad \omega = R/L,
\end{equation*}
with $B_1$ and $B_2$ being dimensionless functions depending also on the approximation scheme  and on the type of boundary conditions (see Fig.~\ref{fig:1}).  Our numerical experiments confirmed these statements {and found no more than one of such sector on the $(h, \tau)$ plane.}

\section{Results of Numerical Experiments}
\label{se.resul}

\subsection{Rod and Scheme Parameters, Initial Conditions, and Reference Solution}
\label{subse.rod}
For our experiments we choose rod parameters that are similar to steel: $\rho= 7860\,\mathrm{kg \: m^{-3}},\;E=210\times 10^9\:\mathrm{Pa}$. The radius of the rod is $R=10^{-3}\:\mathrm{m}$. {Therefore, $C = E\,R^2\,\rho^{-1} \approx 26.717557 \, \mathrm{m}^4\, \mathrm{s}^{-2}$ and $D = R^2 = 10^{-6} \, \mathrm{m}^2$. The length of the rod is $L = 1 \: \mathrm{m}$, and $\omega = R/L = 10^{-3}$}. We set the integration time to $T = 0.3 \: \mathrm{s}$. The point $(h,\,\tau)$ is located in the white areas of Fig.~\ref{fig:1}.

There is an infinite set of such {ADTBCs} gauges that could be obtained by solving Syst.~\eqref{5.4} and Syst.~\eqref{5.5}. The obtained coefficients depend on all parameters of the model. Further we consider several sets of polynomial degrees used in {ADTBCs} construction Eq.~\eqref{5.and}, \eqref{2zv} and \eqref{3zv}. {Therefore,} $h$ and $\tau$ {should be} such that all boundary conditions would exist and not break the stability property.

We choose step $h$ with respect to $x$ to be equal to $0.02\:\mathrm{m}$, and step $\tau$ with respect to $t$ to be equal to $1.6\times 10^{-4}\:\mathrm{s}$ In this case, dimensionless parameters are $\nu \approx 4.274809$ and $\mu = 0.0025$.

\begin{table}[b]
    \centering
    \begin{tabular}{cll}
    \toprule
    \textbf{Parameter} & Value & Dimension \\
    \midrule
      $\rho$ & $7\,860$ & $\text{kg}\,\text{m}^{-3}$ \\
     $E$ & $210\times10^9$ & $\text{Pa}$ \\
     $R$ & $10^{-3}$ & $\text{m}$ \\
     $L$ & $1$ & $\text{m}$ \\
     $h$ & $0.02$ & $\text{m}$ \\
     $\tau$ & $1.6\times10^{-4}$ & $\text{s}$ \\
     $T$ & $0.3$ & $\text{s}$ \\
     $\nu$ & $\approx 4.274809$ & -- \\
     $\mu$ & $0.0025$ & -- \\
    \bottomrule
    \end{tabular}
    \caption{Values of parameters used in numerical simulations.}
    \label{tab:model_params}
\end{table}

\textbf{Note 13.} The physical dimension of the solution is length. However, it may be multiplied by an arbitrary factor, since finite-difference 
Eq.~\eqref{2.2} and {ADTBCs} \eqref{2zv} and \eqref{3zv} are linear. Therefore, the absolute values of the ordinate axis on all figures are optional.

We set initial conditions for Eq.~\eqref{2.1} as 
\begin{equation}
    \label{init_cond_exp}
    u(0, x) = \frac{x}{{\sqrt{\pi\cdot 0.02}}} \, \exp \left(-\frac{x^2}{0.02}\right), \quad \partial_t u(0, x) = 0, \quad x \in \left[-L/2,\, L/2\right].
\end{equation}

Therefore, for Eq.~(\ref{2.2}) we have initial conditions $u_0(x_i) = u(0, x_i)$, and $u_\tau(x_i) = u(\tau, x_i)$ is calculated using the algorithm proposed in \ref{app_init}. The initial condition $u_0(x)$ is close to zero at the ends of the segment $\left[-L/2,\,L/2\right]$.

\begin{figure}[ht]
    \centering
    \includegraphics[scale = 0.4]{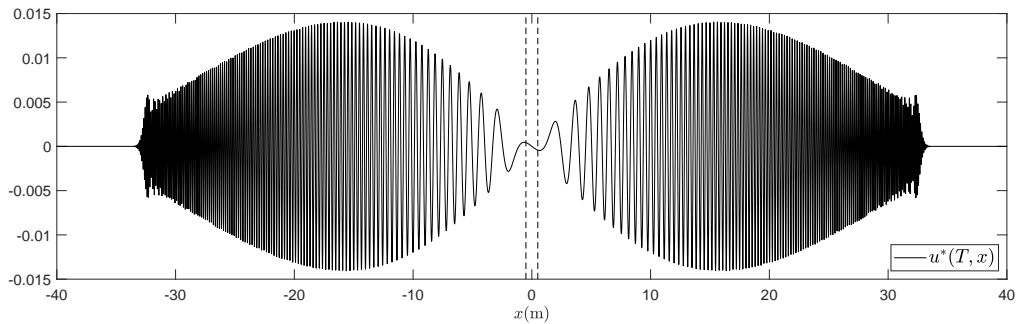}
    \caption{The reference solution $u^*$ on the very  {extended segment $[-40L,\,40L]$ with initial conditions~\eqref{init_cond_exp} at the final time moment $T = 0.3$}. Two vertical dash lines indicate the borders of the considered segment {$x\in[-L/2,\,L/2]$}.}
    \label{fig:2}
\end{figure}

We define a reference solution $u^*(t,\,x)$ on the extended segment $[-40L,\,40L]$, where we use the simplest (Dirichlet + Neumann) homogeneous boundary conditions. The initial rod's perturbation dissipates from the small segment. However, the plot in Fig.~\ref{fig:2} confirms that the boundary conditions cannot significantly influence the solution at the segment $[-0.5L,\,0.5L]$ during the integration period $t \in [0,\, 0.3]$.

The norms of the reference solution $u^*(t,\,x)$ decrease with time $t$ {as $x\in\left[-L/2,\,L/2\right]$}, and the perfect DTBCs must provide the decrease of the corresponding norms of the difference between $u^*$ and the solution of the mixed problem.

The dynamics of obtained solutions and the reference solution, as well as $\mathbf{C}$-norm, are presented in {online version in} \ref{app.anim}.

\subsection{Basic Version of Transparent Boundary Conditions}
\label{subse.basic}

There is an infinite set of choices of the degrees of polynomials $\langle P_1, \, Q_1, \, R_1, \, S_1 \rangle$ and $\langle P_2, \, Q_2, \, R_2, \, S_2 \rangle$, and therefore, an infinite number of corresponding {ADTBCs}. Moreover, it is essential to check the solvability of Syst.~\eqref{5.4} and Syst.~\eqref{5.5}. If at least one of the systems does not have a unique solution for some set of polynomials degrees, then  we cannot construct the {ADTBCs} for this set of degrees.

{When the rational approximation in Syst.~\eqref{5.and} is performed, the stability of the Crank -- Nicolson scheme might be lost. Therefore, one should account for solvability of Syst.~\eqref{5.4} and Syst.~\eqref{5.5}, and check that the stability of the mixed problem is preserved. If the {ADTBCs} exist for the  specific model's parameters and polynomial degrees, we apply them in computations and check if the result remains stable. Rational approximations~\eqref{5.4} and \eqref{5.5} with finite polynomial degrees do not guarantee obtaining completely transparent boundaries. To determine if {ADTBCs} are reasonable, we investigate the error of the obtained solution $u$ by comparing it to the reference solution $u^*$.

Syst.~\eqref{5.4} for $k=1$, and Syst.~\eqref{5.5} for $k=2$ {consist of two similar} equations for real and imaginary parts with the smallness order $K$. We also have two normalisation conditions on coefficients. Thus, the total number of unknown coefficients of these four polynomials that we have determined from Syst.~\eqref{5.4}, as well as Syst.~\eqref{5.5}, should be even.

As an example, let us consider the set of polynomial degrees
\begin{equation}
\label{7.1.1}
    \text{deg}\,P_k = \text{deg}\,Q_k = 4, \; \text{deg}\,R_k = \text{deg}\,S_k = 8, \quad k=1,\,2.
\end{equation}

{Here and further we consider equal sets of polynomial degrees of boundary conditions for both border ($k=1$) and preborder ($k=2$) points}. The coefficients that correspond to the solutions of Syst.~\eqref{5.4} and Syst.~\eqref{5.5} for degree sets (\ref{7.1.1}) are presented in Table~\ref{tab:2}.

\begin{table}[ht]
    \centering
    \begin{tabular}{l|rrrr | rrrr}
         & $P_1$ & $Q_1$ & $R_1$ & $S_1$ & $P_2$ & $Q_2$ & $R_2$ & $S_2$  \\
        \toprule
        1&   1&   0&    -0.555979&   0,278657&   0&   1&    -0.925737&   0.301010 \\
        $\omega$&   -1.039354&  -1.064260&  0.925512&   -0.300505&  -0.039239&  -1.498177&  0.962232&   -0.272787   \\
        $\omega^2$& 1.040798&   0.175892&   -0.343658&  0.205584&   -0.057023&  1.346122&   -0.918314&  0.289728    \\
        $\omega^3$& -0.484423&  -0.688193&  1.007943&   -0.361839&  0.240692&   -1.187154&  0.993006&   -0.295379   \\
        $\omega^4$& 0.217631&   -0.187829&  0.258996&   -0.095354&  -0.007746&	0.054903&   0.027530&   -0.020261   \\
        $\omega^5$&         &           &   0.101158&   -0.063710&           &           &   0.039188&	-0.023854   \\
        $\omega^6$&         &           &	0.008250&   -0.016540&           &           &	0.004642&	-0.006821   \\
        $\omega^7$&	        &           &	-0.014938&  0.002764&            &           &   -0.005037&  0.000709   \\
        $\omega^8$&         &           &   -0.005839&  0.002373&            &           &	-0.002124&	0.000827    \\
    \end{tabular}
\caption{The coefficients of the {ADTBCs} are obtained from Eqs.~(\ref{2zv}) and (\ref{3zv}) for the sets of polynomial degrees $\langle P_k, \, Q_k, \, R_k, \, S_k \rangle = \langle 4, \, 4, \, 8, \, 8 \rangle$ at $k = 1,\,2$.}
\label{tab:2}
\end{table}

\textbf{Note 14.} We provide values in the tables up to 6 decimal places. Our numerical experiments showed that by using the ADTBC coefficients without the 6-th decimal place the error increases just slightly, whereas without the 5-th decimal place the increase is significant.

To evaluate the dynamics of the  error of the obtained solution $u$ of Eq.~(\ref{2.2}) with {ADTBCs}, we use the common logarithm (base $10$) of three norms:
\begin{enumerate}[{\bf a)}]
    \item common logarithm of the {Hamiltonian approximation $\hat{\mathcal{H}}$} of solutions' difference:
    
    \begin{equation*}
        \log_{10} \sqrt{\hat{\mathcal{H}} \left[u(t, x) - u^*(t, x)\right]},
    \end{equation*}

    \item common logarithm of the Chebyshev norm $\mathrm{C}[-L/2,\,L/2]$ of solutions' difference, i.e.
    \begin{equation*}
        \log_{10} \left[\max_{j=0,1,\ldots,N} \left| u(t, x_j) - u^*(t, x_j) \right|\right],
    \end{equation*}
    
    \item common logarithm of $\mathbf{L}^2$-norm of solution's deference
    \begin{equation*}
        \log_{10}\left[ \|u(t, x) - u^*(t, x)\|_2\right],
    \end{equation*}
\end{enumerate}
where $\|u(t, x) - u^*(t, x)\|_2^2 = \int_{-L/2}^{L/2} \left(u(t, x) - u^*(t, x)\right)^2 \, \mathrm{d}x$, which is approximated by the standard trapezoidal method.

The results of our numeric experiments  with different {ADTBCs} sets are presented in Fig.~\ref{fig:3}. 

\begin{figure}[ht]
    \centering
    \begin{subfigure}[b]{.33\linewidth}
        \includegraphics[width = \textwidth]{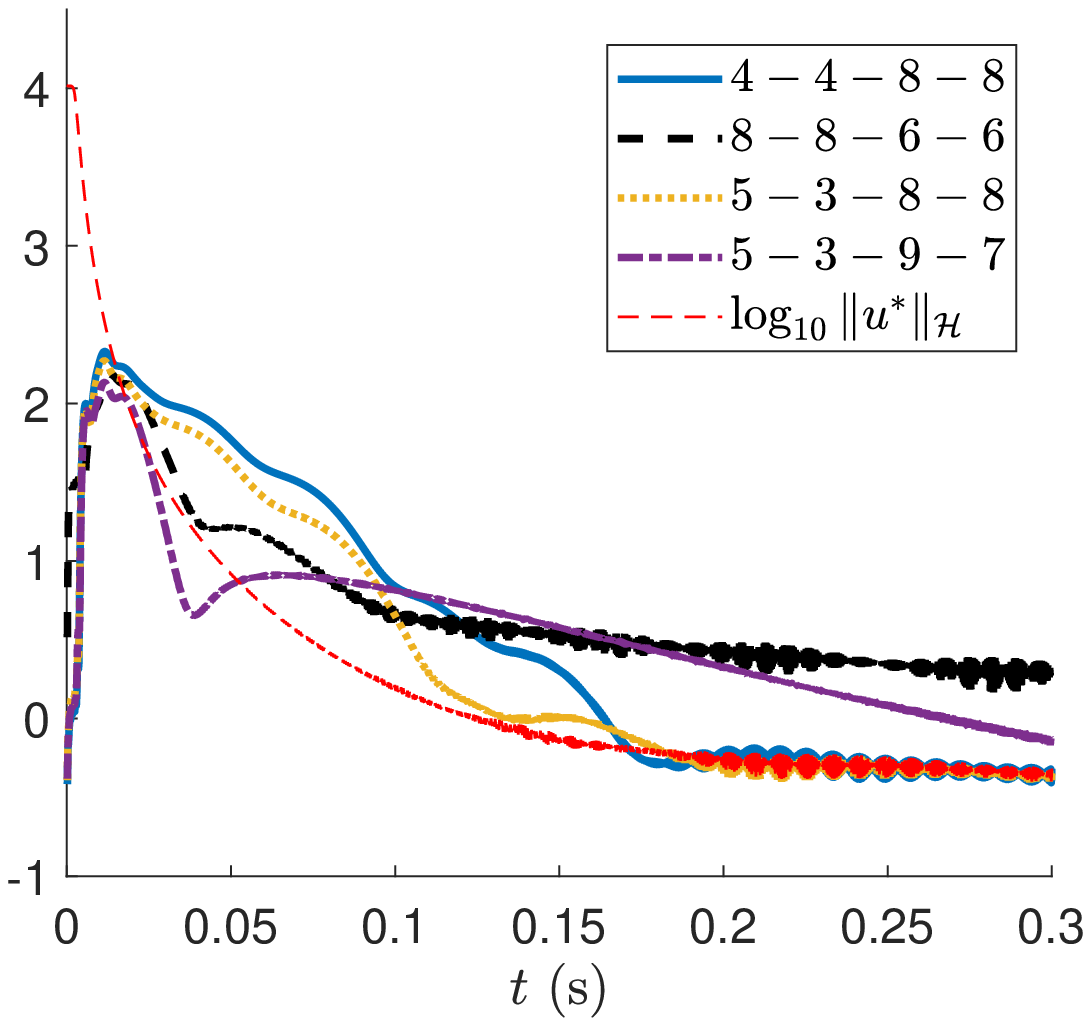}
        \caption{\label{fig:3a}}
    \end{subfigure}
    \begin{subfigure}[b]{.33\linewidth}
        \includegraphics[width = \textwidth]{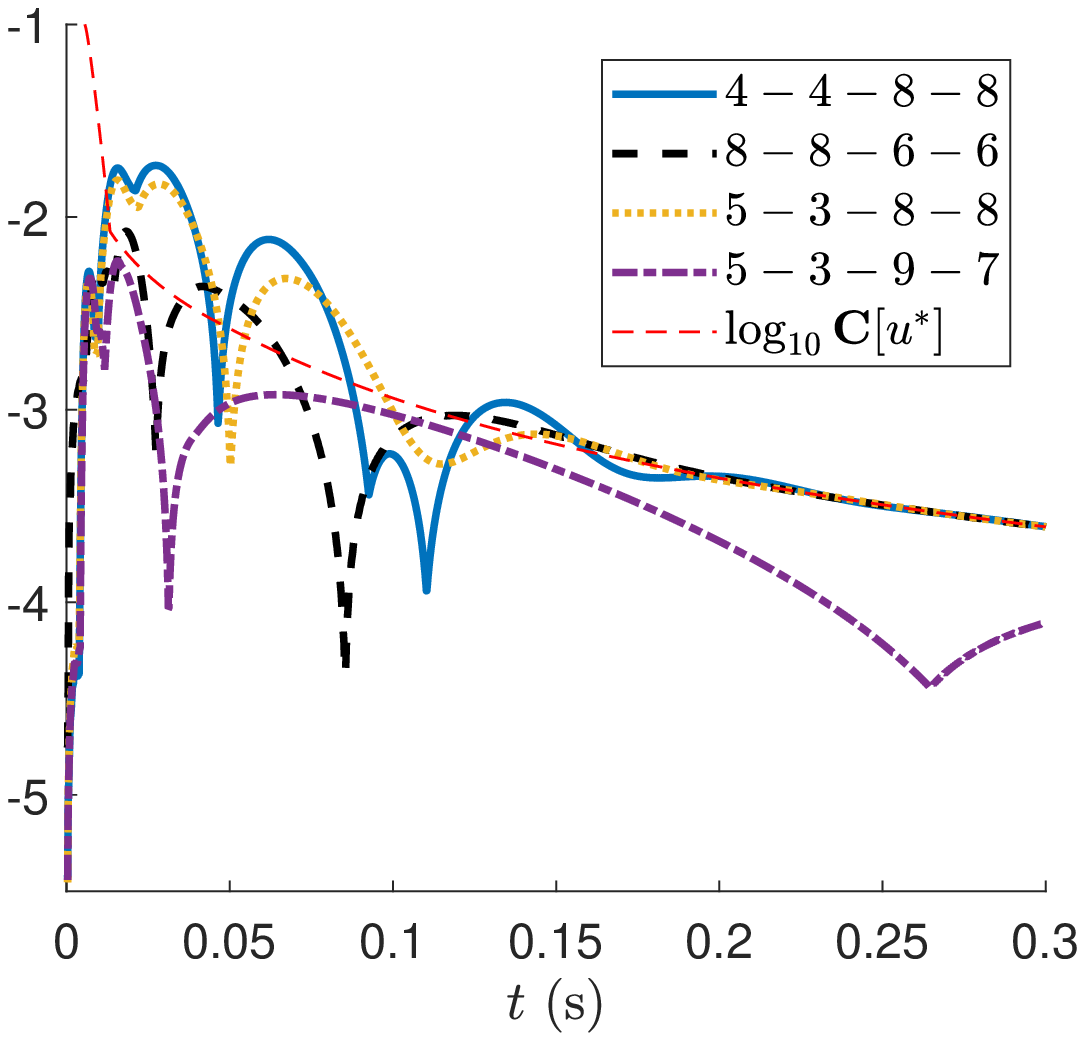}
        \caption{\label{fig:3b}}
    \end{subfigure}
    \begin{subfigure}[b]{.33\linewidth}
        \includegraphics[width = \textwidth]{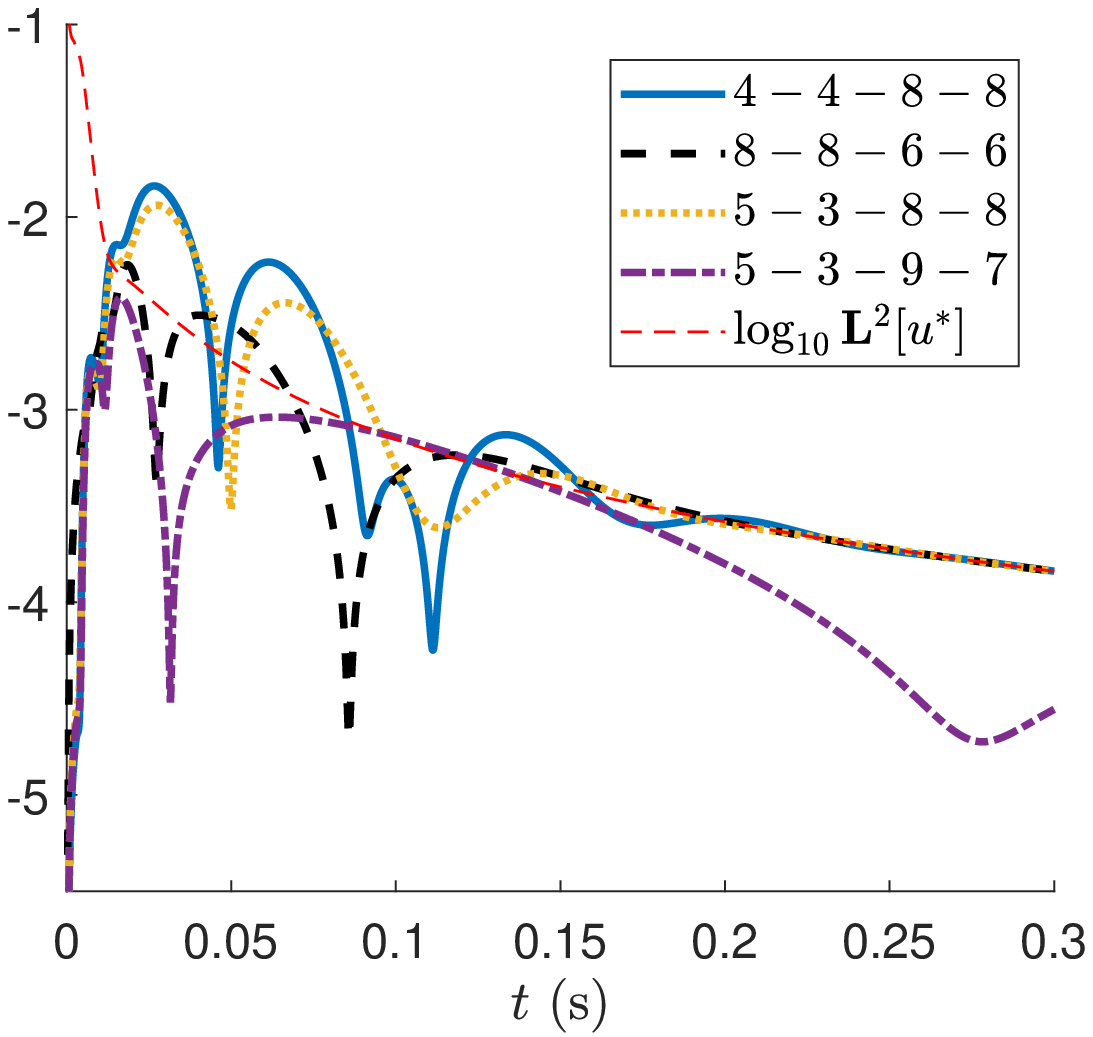}
        \caption{\label{fig:3c}}
    \end{subfigure}
    \caption{{Solid, thick dashed, dotted and dash-dotted lines represent common logarithm of (a) $\hat{\mathcal{H}}$, (b) \textbf{C}-norm, (c) $\mathbf{L}^2$-norm of the difference between the reference solution $u^*$ and solutions with {ADTBCs}. The narrow dashed line is common logarithm of (a) $\hat{\mathcal{H}}$, (b) \textbf{C}-norm, (c) $\mathbf{L}^2$-norm of the reference solution $u^*$. Initial conditions are defined by formula \eqref{init_cond_exp}.}}
    \label{fig:3}
\end{figure}

From Fig.~\ref{fig:3a} it follows that a part of the energy of the {reference} solution decreases approximately 20 000 times compared to the initial moment, because waves go away from the small segment. However, $\hat{H} (t)$ decreases slowly and non-monotonically  when $t\to \infty$, see Fig.~\ref{fig:4a}, where the  zoomed in fragment shows the dynamics of $\hat{H} (t)$ at $t\approx 0.25$. The oscillation of this functional occurs over time, the difference between two local maxima on the plot is equal to $2\tau$. The amplitude of the oscillations decreases very slowly, see Fig.~\ref{fig:4b}.

\textbf{Note 15.} Let us assume that at large times $t$, the oscillation energy is distributed almost uniformly over a segment and expands at about a constant speed on the straight line $x\in\mathbb{R}$. Then the energy that is concentrated on the segment $[- L/2, L/2]$ should decrease approximately as ${\bf O}\left(t^{-1}\right)$, and the energy norm of the solution as ${\bf O}\left(t^{-1/2}\right)$. Our evaluation shows that for this initial condition, the energy norm at $t\to \infty$ is estimated asymptotically as $0.26\cdot t^{-1/2}$, see Fig.~\ref{fig:4b}.}

However, the {ADTBCs} with polynomial degrees $\langle 4, 4, 8, 8 \rangle$, as shown on Fig.~\ref{fig:4}, significantly reduce the part of the energy on the segment $[-L/2,\,L/2]$ in comparison to the reference solution. The difference of these energy norms is so significant that we needed to multiply the {ADTBCs} solution by 10 to place it in Fig.~\ref{fig:4a}.

{\textbf{Note 16.} Also, we found that if the initial function $u_0 \equiv u(0, x)$ satisfy the Eq.~\eqref{eq:zeroth_int}, then we can set additional requirements on the coefficients in Syst.~\eqref{5.4} and \eqref{5.5}, which may reduce the error even further. However, the motivation behind this modification of DTBCs and ADTBCs does not have a physical interpretation, and the resulting error is not necessarily lower than in the standard method (described in Subsect.~\ref{subse.rational}). We present the modification and results in \ref{app.modif}.}

\begin{figure}[ht]
    \centering
    \begin{subfigure}[b]{.45\linewidth}
        \includegraphics[width = \textwidth]{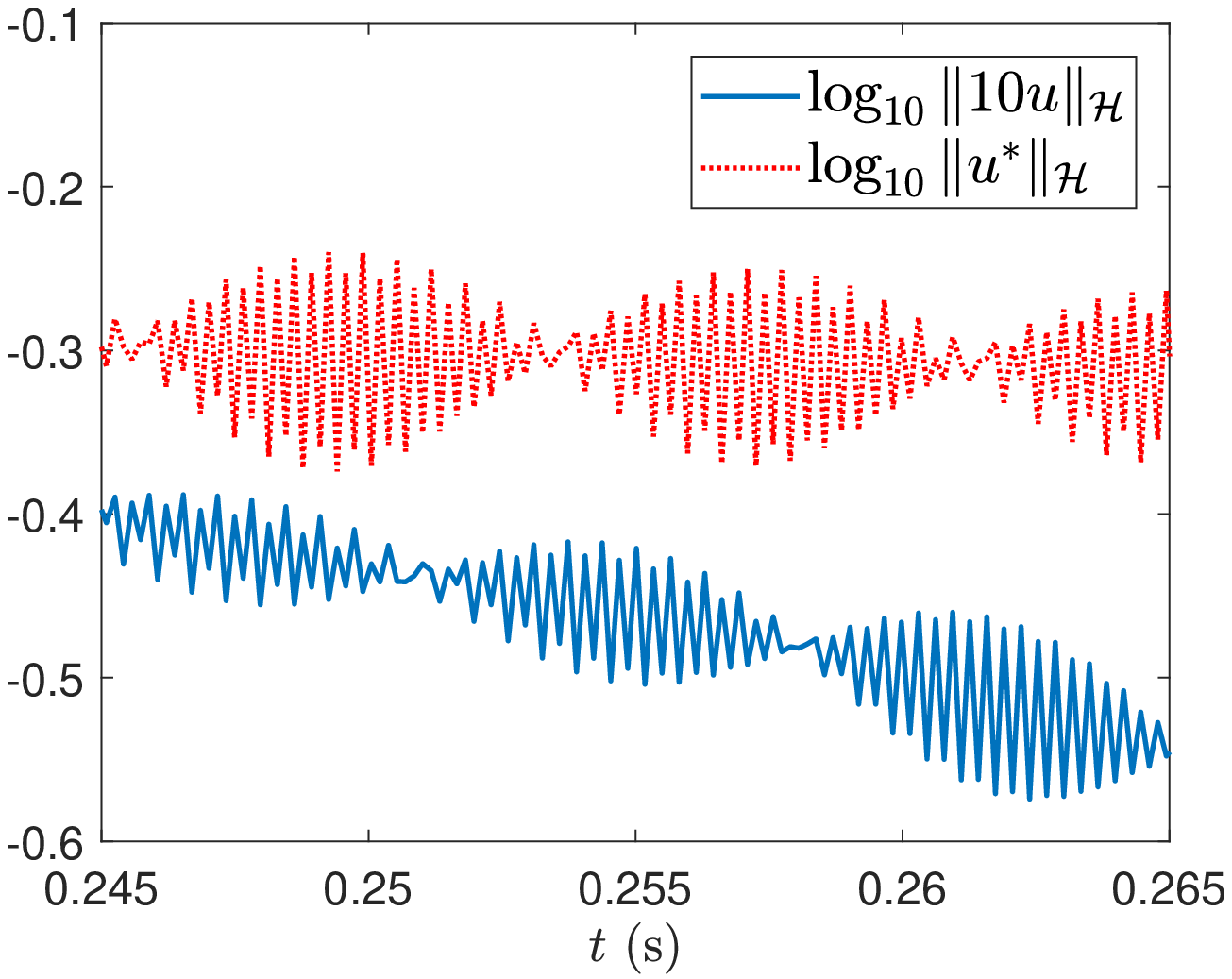}
        \caption{\label{fig:4a}}
    \end{subfigure}
    \begin{subfigure}[b]{.45\linewidth}
        \includegraphics[width = \textwidth]{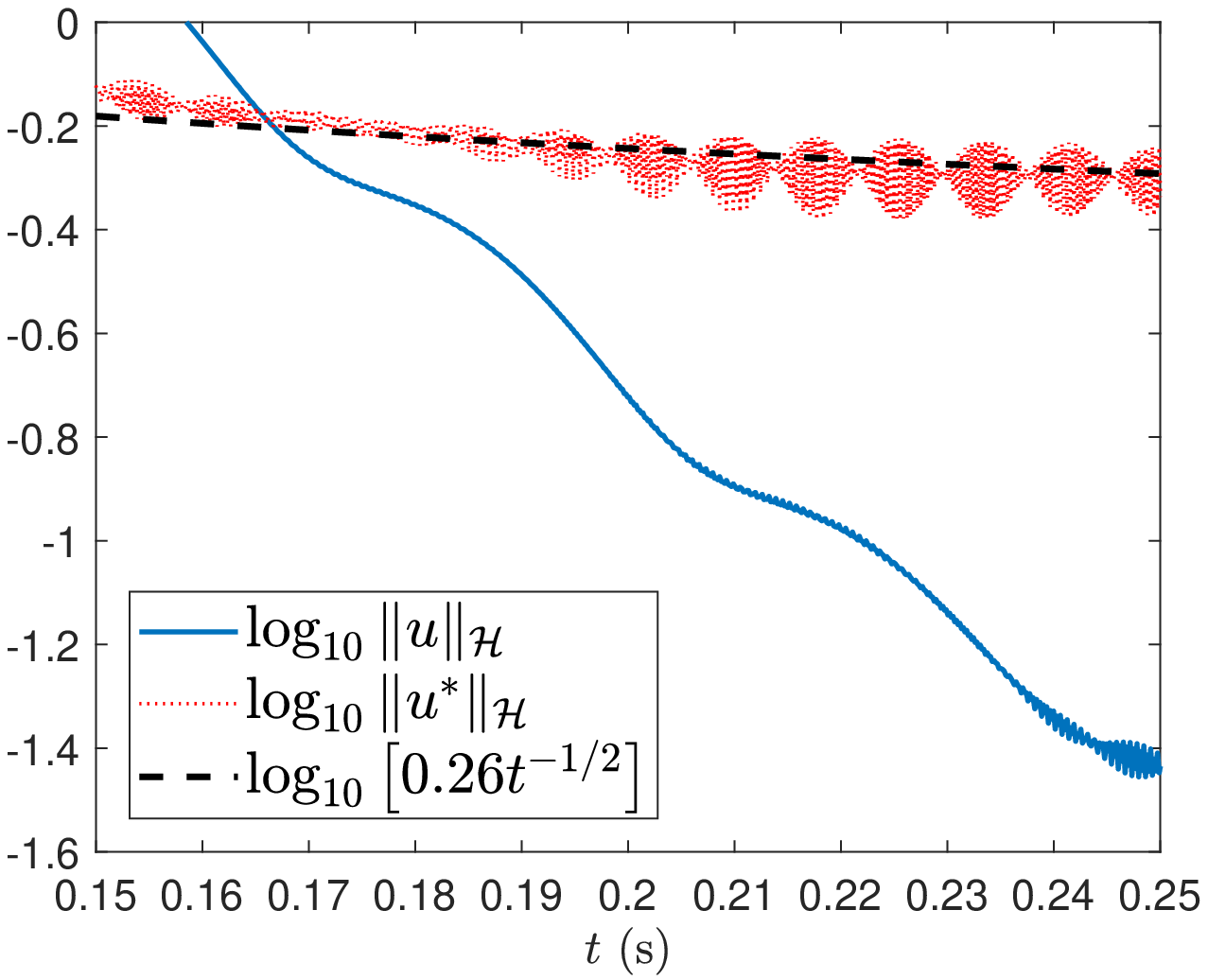}
        \caption{\label{fig:4b}}
    \end{subfigure}
    \caption{{(a) Common logarithm of {$\hat{\mathcal{H}}$} of (solid line) the solution (multiplied by $10$) obtained with {ADTBCs} with polynomial degrees $\langle4,4,8,8\rangle$ and (dot line) the reference solution $u^*$ on time interval $t \in [0.245,\,0.265]$. (b) The decrease of energy in the period when most of it left the small segment, but has not yet reached the edges of the large segment. The dash line shows the asymptotic of the reference solution. Initial conditions are defined by formula~\eqref{init_cond_exp}.}}
    \label{fig:4}
\end{figure}

\subsection{Comparison of Transparent Boundary Conditions with Various Versions of `Usual' Homogeneous Ones}
\label{subse.comparison}
In practice, simple homogeneous boundary conditions (i.e. Dirichlet and Neumann) are usually used when there is no information about physical processes on the border. These conditions lead to the partial or complete reflection of outgoing waves, back into {the computational domain} (sometimes with increased amplitude). On the contrary, {ADTBCs} that are calculated using our vectorial rational approximation techniques have a low reflection level.

\begin{figure}[ht]
    \centering
    \begin{subfigure}[b]{.33\linewidth}
        \includegraphics[width = \linewidth]{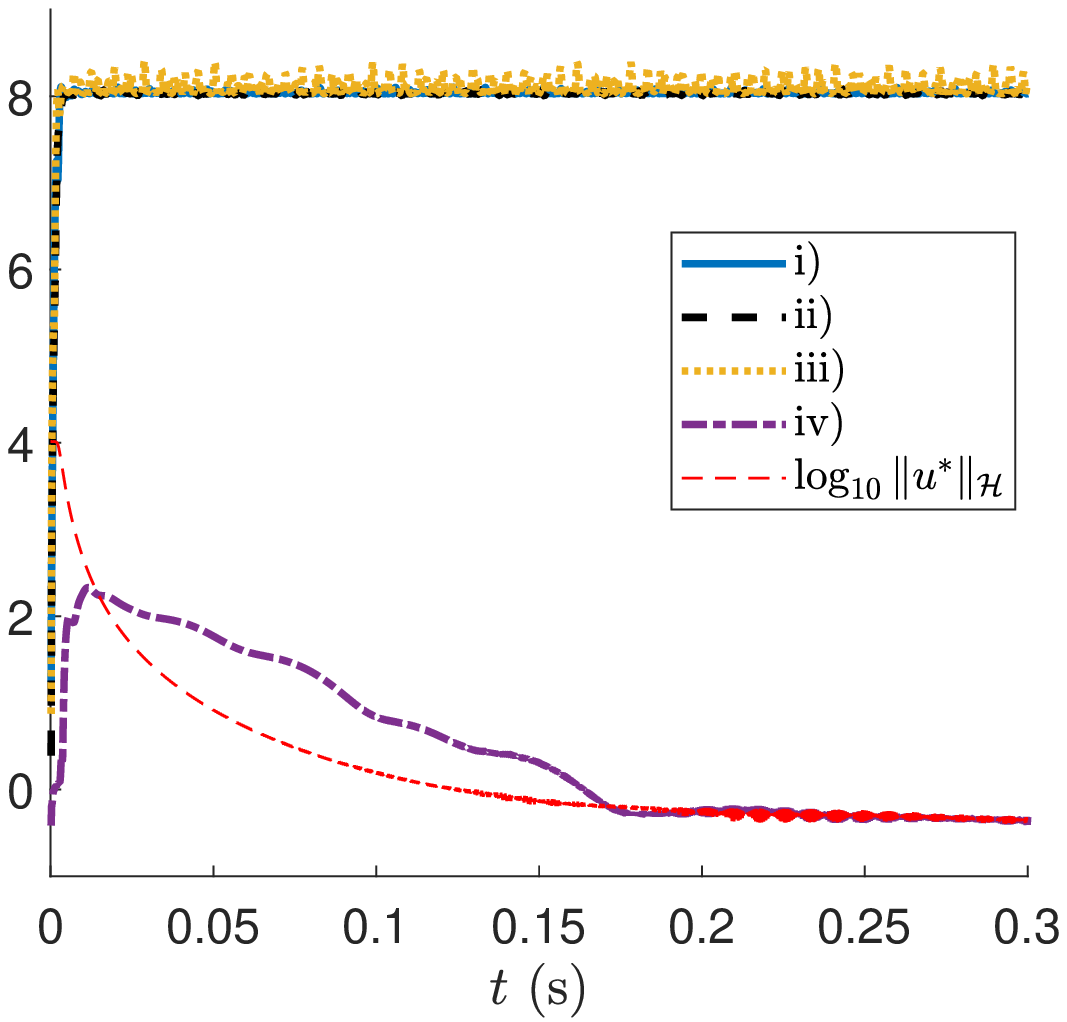}
        \caption{\label{fig:7a}}
    \end{subfigure}
    \begin{subfigure}[b]{.33\linewidth}
        \includegraphics[width = \linewidth]{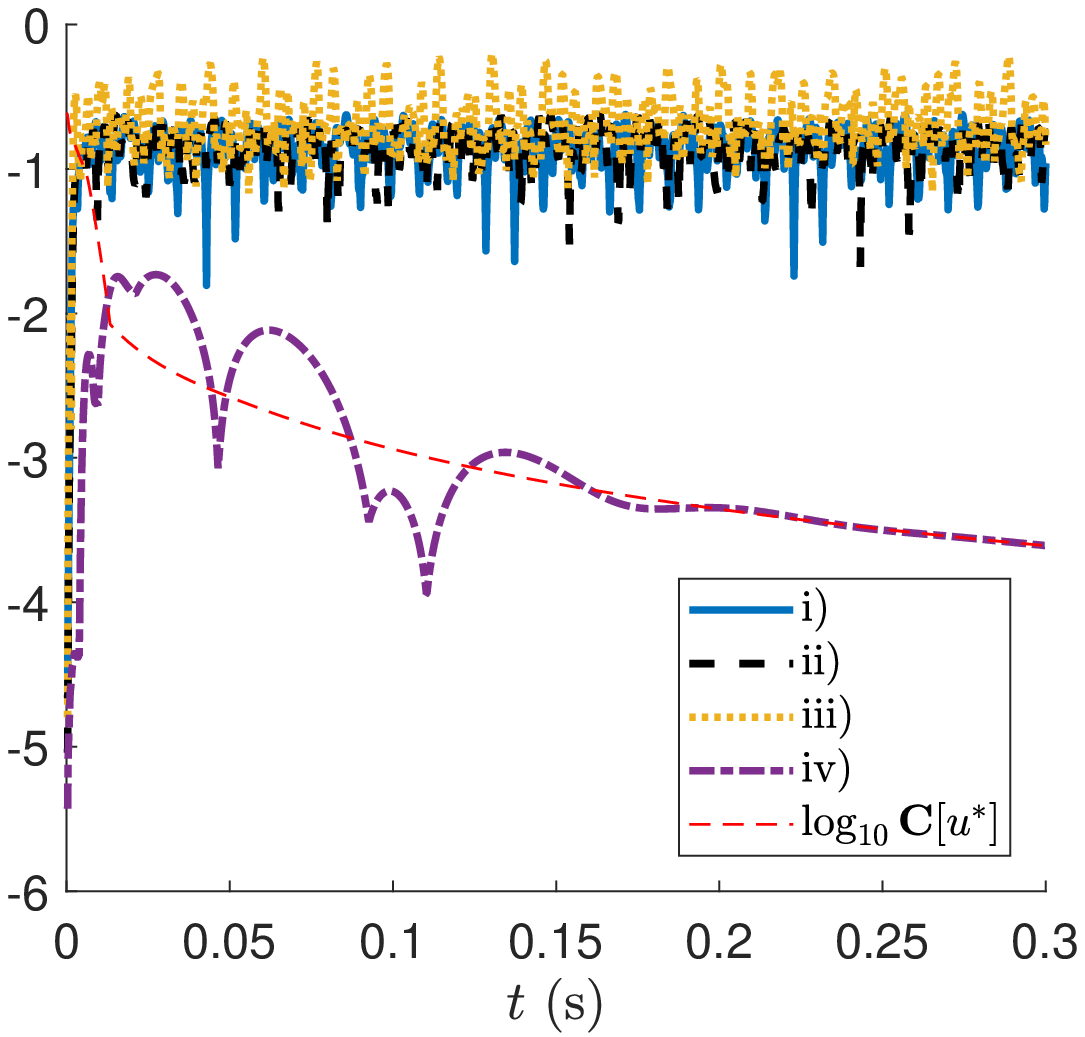}
        \caption{\label{fig:7b}}
    \end{subfigure}
    \begin{subfigure}[b]{.33\linewidth}
        \includegraphics[width = \linewidth]{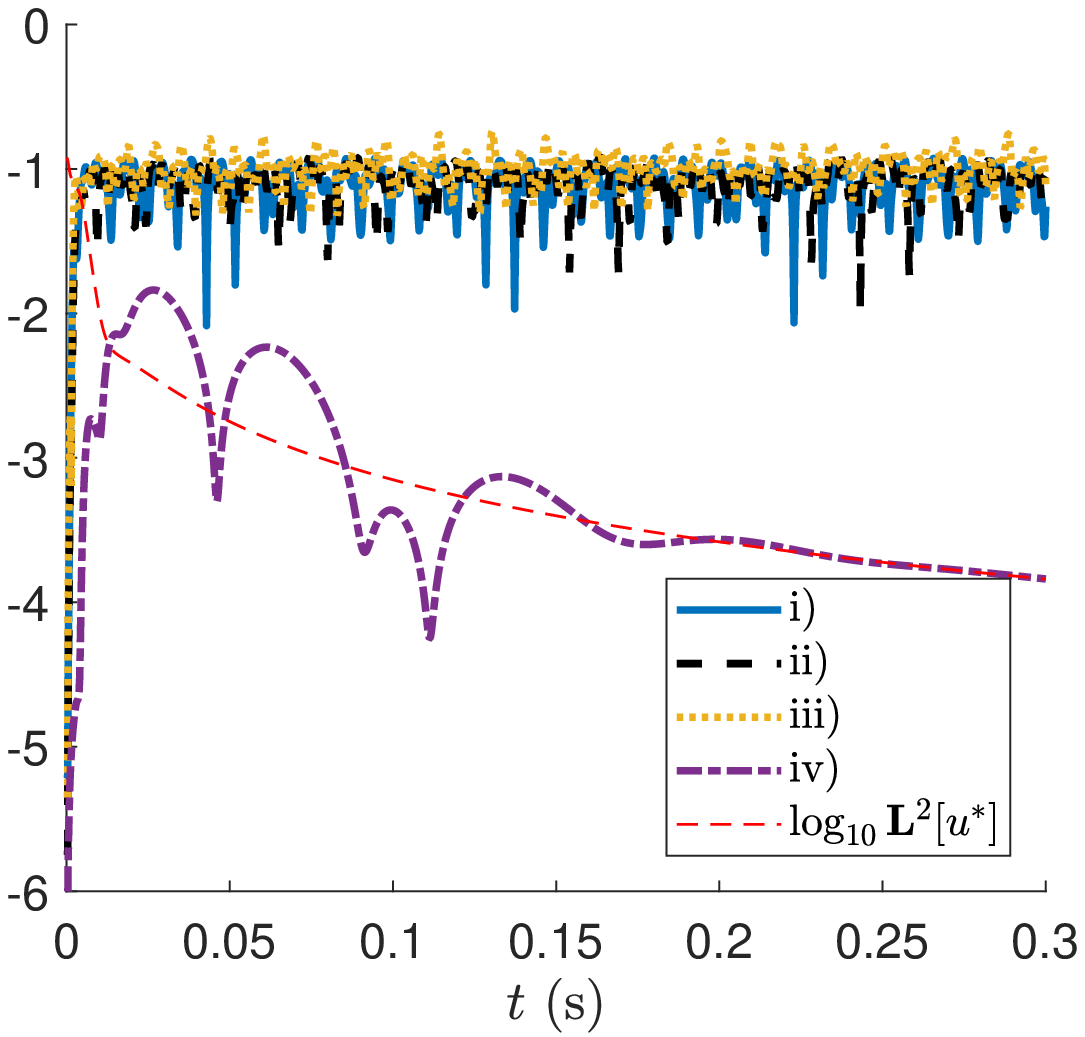}
        \caption{\label{fig:7c}}
    \end{subfigure}
    \caption{The dashed, dotted and dash-dotted lines represent {common logarithm of (a) {$\hat{\mathcal{H}}$}, (b) \textbf{C}-norm  and (c) $\mathbf{L}^2$-norm of the difference between reference solution $u^*$ and solution  with `usual' boundary conditions. The narrow dash line is common logarithm of (a) {$\hat{\mathcal{H}}$}, (b) \textbf{C}-norm  and (c) $\mathbf{L}^2$-norm of the reference solution $u^*$. Initial conditions are defined by formula~\eqref{init_cond_exp}.}}
    \label{fig:7}
\end{figure}

Fig.~\ref{fig:7} shows the dynamics of solutions' errors that are calculated using various `usual' boundary conditions:
\begin{enumerate}[i)]
    \item $u|_{\Gamma} = 0, \, \left. \frac{\partial u}{\partial x} \right|_\Gamma = 0 \implies u_0^n = u_1^n = 0$,
    \item $u|_{\Gamma} = 0, \, \left. \frac{\partial^2 u}{\partial x^2} \right|_\Gamma = 0 \implies u_0^n = 0, \, u_1^n = u_2^n / 2$,
    \item $\left. \frac{\partial^2 u}{\partial x^2}\right|_{\Gamma} = 0, \, \left. \frac{\partial^3 u}{\partial x^3} \right|_\Gamma = 0 \implies u_0^n = 3u_2^n - 2u_3^n, \, u_1^n = 2u_2^n - u_3^n$,
    \item {ADTBCs} with polynomial degrees $\langle P_k,\, Q_k,\, R_k,\, S_k \rangle = \langle 4,\, 4,\, 8,\, 8 \rangle$, $k = 1,\,2$.
\end{enumerate}

All these homogeneous boundary conditions i) - iii) lead to a significant reflection of outgoing waves back into the computational area, see Fig.~\ref{fig:7}. The obtained solutions almost immediately differ from the reference solution, whereas the solution obtained with {ADTBCs} (iv) stays close to $u^*$ during the integration time.

\begin{figure}[h]
    \centering
    \includegraphics[scale = 0.5]{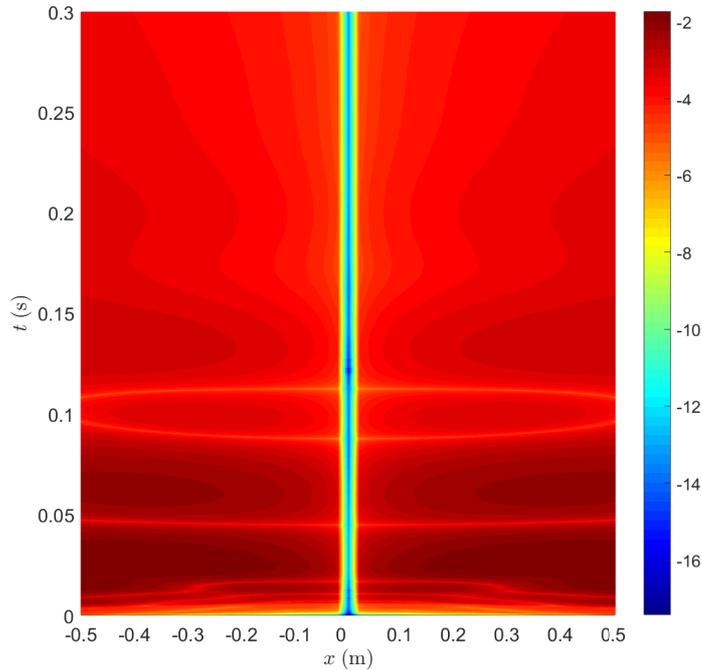}
    \caption{{Common logarithm of absolute difference} between reference solution $u^*$ and the mixed problem solution with {ADTBCs} with polynomial degrees $\langle P_k,\, Q_k,\, R_k,\, S_k \rangle = \langle 4,\, 4,\, 8,\, 8 \rangle$, $k = 1,\,2$ at different time moments $t$. {Initial conditions are defined by formula~\eqref{init_cond_exp}.}}
    \label{fig:8}
\end{figure}

The smoothed evolution of the solution with {ADTBCs} (with degrees $\langle 4, 4, 8, 8 \rangle$) is shown in Fig.~\ref{fig:8}. Lines where $|u(t,x)-u^*(t,x)|=0$ are generally not visible, except in the middle of the segment, where both solutions are small, because they are odd for all $t$. Fig.~\ref{fig:8} shows that for large $t$, the difference $|u(t,x)-u^*(t,x)|$ is slowly decreasing. According to Fig.~\ref{fig:4}, it is because the reference solution slowly decreases, but the solution with {ADTBCs} at this time is {orders of magnitude less}.

\section{Discussion}
\label{se.discussion}
{ADTBCs} might be used in the mathematical modelling of processes on a {finite} area, when it is certain that external processes do not have any essential impact on the interior. On the other hand, {ADTBCs} could suspend all fluctuations in the finite area without using any high computational  viscosity on some borders' vicinity.

The problem of constructing transparent conditions becomes more complex if the coefficients of the partial differential (or finite-difference) equation are variable, and/or the bounded computational domain $V$ has a sophisticated boundary. In such cases, it is usually not possible to apply the method of separation of independent variables. {In \cite{ZAITSEV2007}, for the anisotropic elasticity, it is proposed to generalise the approach \cite{Sofronov1992} of fast analytical TBC for the wave equation to the case of variable coefficients. The quasi-analytical ADTBC operator \cite{ZAITSEV2007} is generated numerically and its accuracy depends on the number of basis functions representing the solution at the open boundary. The problem of arbitrary computational domain can be solved by immersion into a sphere with ADTBC using domain decomposition conception as proposed in \cite{Sofronov2007}.} In the works \cite{engquist1977absorbing, engquist1979radiation, gordin1979diss, gordin2000mathematical, gordin1987mathematicalb,  gordin1982application,  gordin1979upper, gordin1978projectors, gordin1978mixed, gordin1977}  (and also in this article), only problems where the variable separation is possible were considered. On the other hand, here the construction of DTBCs and ADTBCs for the difference problem is performed directly and not by means of a difference approximation of the boundary operators for the differential problem.

{In our work we propose the algorithm of determining {ADTBCs} operators, which can be done for fixed parameters of the rod and time-space steps. However, if the parameters or pair $(h, \tau)$ are changed, a numerical recalculation is needed. For the Schr\"odinger equation to recalculate DTBC operators for other pairs $(h, \tau)$ a transformation rule was proposed by using one and the same rational approximation, see \cite{arnold2003}.}

A characteristic feature of DTBCs is non-locality with respect to time. Values of the solution in the vicinity of a border at previous time steps are required. {We used a vectorial rational approximation generalising the Hermite -- Pad\'e approximation. It allows us to reduce this number of time steps in the convolutions}. {In spatially multidimensional  finite-difference models, where the corresponding DTBCs are non-local with respect to variables that are tangential to the area's boundary, vectorial rational approximation can also be applied, which would significantly reduce computational costs,} see \cite{engquist1977absorbing, gordin1977, engquist1979radiation, gordin1978mixed, gordin1978projectors, gordin1979upper, gordin1982boundary, gordin1987mathematical, gordin1987mathematicalb, gordin2000mathematical, gordin2010mathematics, arnold2005approximation, arnold2003, ZAITSEV2007, Sofronov1992}.

\section{Conclusion}
\label{se.con}
In this paper, we have constructed {Approximate Discrete Transparent Boundary Conditions (ADTBCs)} for the implicit finite-difference scheme that approximates the fourth order differential equation with respect to space. Both equations (differential and finite-difference) require two boundary conditions on each end. The considered differential equation is more sophisticated {than} many classic mathematical physics equations, because it is not resolved with respect to the highest derivative of a solution with respect to time (i.e. it does not belong to Cauchy -- Kovalevskaya type).

Here, the {ADTBCs} were constructed for the finite-difference Crank -- Nicolson implicit scheme for the transverse vibrations  equation of a rod (beam) with a circular cross section. {ADTBCs} provide a solution of a mixed initial-boundary value problem on a segment that is close to the solution on the infinite domain. 

{We proved the absolute stability of the Crank -- Nicholson scheme for the Cauchy problem, and experimentally verified the conditional stability of the mixed boundary value problem with {ADTBCs}.} It is shown that the stability regions depend on the rational approximation (i.e. sets of polynomial degrees), and are bound by two parabolas on the $(h,\tau)$ plane.

It is shown that `usual' {(Dirichlet and  Neumann)} homogeneous boundary conditions do not have this `transparency' property. The need of such {ADTBCs} is seen in many scientific and technical applications. The approach may be applied, when we need to imitate a non-zeroth background solution (at a large area) and a forcing $f$. In these cases {ADTBCs} will be non-homogeneous.

All {DTBCs and ADTBCs} derived here, are characterised by lots of numerical parameters (physical parameters of a rod, degrees of approximating polynomials in Syst.~\eqref{5.and}, space $h$ and time $\tau$ steps). Resulting {ADTBCs} should be defined with at least five decimal places. We describe the algorithm of parameter determination (symbolic computations were used), which is the main result of this paper.

The proposed algorithms of {ADTBCs} construction can be used for various evolutionary linear equations or systems and their finite-difference approximations. However, recalculations of all {coefficients and formulae derivations for} the {ADTBCs} are required, if other finite-difference approximation schemes are used.

We also present the algorithm (based on the compact finite-difference scheme) of the initial functions calculation that provides a high-order of approximation with respect to time, for the implicit finite-difference equation.

\section*{Acknowledgements}

We are cordially grateful to Ph. L. Bykov for useful discussions in the course of our work. We are very grateful to anonymous reviewers, whose comments helped us to correct inaccuracies, confusions and deficiencies.

The article was prepared within the framework of the Academic Fund Program at the National Research University Higher School of Economics (HSE University) in 2018 -- 2019 (grant \textnumero~18-05-0011), in 2020 -- 2021 (grant \textnumero~20-04-021)  and by the Russian Academic Excellence Project ``5-100''.

\appendix

\section{{Expansion of Functions $\eta_{1,2}(\omega)$ into Taylor Series}}
\label{app.eta12}
We represent the radicand in Eq.~(\ref{4.6}) with a help of dimensionless parameters $\nu$ and $\mu$ and simplify it:

\begin{equation*}
    \eta_1 (\omega) = \frac{-\beta\left(1 + \omega^2 \right) - \gamma \omega - 
    \sqrt[+]{\left(1 - \omega \right)^2 \left[\left(\mu^2 - 2\nu\right) \omega^2 - 2\mu^2 \omega + \mu^2 - 2\nu\right]}}
    {2\sigma \left(1 + \omega{^2} \right)}.
\end{equation*}

If {$|\omega|< 1$}, we can take out the multiplier from the quadratic root:

\begin{equation*}
    \eta_1 (\omega) = \frac{-\beta\left(1 + \omega^2 \right) - \gamma \omega - 
    \sqrt{\mu^2 - 2\nu} \left(1 - \omega \right) \sqrt[+]{\omega^2 - 2\frac{\mu^2}{\mu^2 - 2\nu} \omega + 1}}{2\sigma \left(1 + \omega^2 \right)}.
\end{equation*}

Let us apply the formula for the  generating function of the Legendre polynomials (see e.g. \cite{abramowitz1972handbook, szeg1939orthogonal, OLVER1974, gordin1978projectors, gordin1982boundary, gordin1987mathematicalb, gordin2000mathematical, gordin2010mathematics}):
\begin{equation*}
\left(\omega^2-2\varepsilon\omega +1\right)^{-1/2} = \sum_{n=0}^\infty \mathrm{P}_n(\varepsilon) \omega^n,
\end{equation*}
where $\mathrm{P}_n(\varepsilon)$ is the Legendre polynomial of degree $n$ at point $\varepsilon$. Here $\varepsilon=\frac{\mu^2}{\mu^2-2\nu},\;|\varepsilon|<1$. We obtain
\begin{equation*}
    \eta_1(\omega) = \frac{1}{2\sigma(1+\omega^2)} \left[ -\beta \left(1+\omega^2\right) - \gamma\omega - \sqrt{\mu^2-2\nu} \left(1-\omega\right) \left(\omega^2 - 2\frac{\mu^2}{\mu^2-2\nu} \omega + 1\right) \sum_{n=0}^\infty \mathrm{P}_n\left( \frac{\mu^2}{\mu^2-2\nu}\right) \omega^n \right].
\end{equation*}

Then we use the formula for geometric progression:
\[
\frac{1}{1+a}=\sum\limits_{k=0}^\infty (-1)^k a^k,
\]
where $a=\omega^2$, and express values $\beta,\,\gamma,\,\delta$ across $\mu$ and $\nu$. We obtain:
\begin{equation}
    \eta_1(\omega) = \frac{1}{\nu} \sum_{k=0}^\infty (-1)^k \omega^{2k} \cdot\left[ \left(\mu+2\nu\right) \left(1+\omega^2\right) - 2\mu\omega - \sqrt{\mu^2-2\nu} \left(1-\omega\right) \left(\omega^2 - 2\frac{\mu^2}{\mu^2-2\nu} \omega + 1\right) \sum_{n=0}^\infty \mathrm{P}_n\left( \frac{\mu^2}{\mu^2-2\nu}\right) \omega^n \right].
\end{equation}

In the same way we obtain
\begin{equation}
\begin{aligned}
    \eta_2(\omega) = \frac{1}{\nu} \sum_{k=0}^\infty (-1)^k \omega^{2k} \cdot\left[ \left(\mu+2\nu\right) \left(1+\omega^2\right) - 2\mu\omega + \sqrt{\mu^2-2\nu} \left(1-\omega\right) \left(\omega^2 - 2\frac{\mu^2}{\mu^2-2\nu} \omega + 1\right) \sum_{n=0}^\infty \mathrm{P}_n\left( \frac{\mu^2}{\mu^2-2\nu}\right) \omega^n \right].
\end{aligned}    
\end{equation}

\section{{Expansion of Functions $\lambda_{i}(\omega)$, $i = 1, \ldots, 4$ into Taylor Series}}
\label{app.lambda}

We can rewrite the functions $\lambda_{1,3}(\omega)$ in the form
\begin{equation*}
    \lambda_{1,3}(\omega)=  \frac{\eta_1(\omega)}{2}\mp \sqrt{\frac{1}{4} \left(\vartheta_1 + r_1(\omega)\right)^2 - 1},
\end{equation*}
and factor the radicand:
\begin{equation*}
    \lambda_{1,3}(\omega)=  \frac{\eta_1(\omega)}{2}\mp \sqrt{\left(\frac{\vartheta_1 + r_1(\omega)}{2} + 1\right) \cdot \left(\frac{\vartheta_1 + r_1(\omega)}{2} - 1\right)}.
\end{equation*}

Then we represent the radicand as a product
\begin{equation*}
    \lambda_{1,3}(\omega) = \frac{\eta_1(\omega)}{2} \mp \sqrt{\frac{\theta_1^2}{4} - 1} \cdot \sqrt{1 + \frac{r_1(\omega)}{\vartheta_1+2}} \cdot \sqrt{1 + \frac{r_1(\omega)}{\vartheta_1-2}},
\end{equation*}
and use the formula for the Taylor series expansion of a square root:
\begin{equation*}
    \sqrt{1+x} = \sum_{n=0}^\infty \frac{(-1)^n \, (2n)!}{(1-2n) \, n! \, 4^n} \, x^n,\quad |x|<1,
\end{equation*}
to obtain the Taylor series for the  characteristic roots in the vicinity of $\omega = 0$
\begin{equation}
    \lambda_{1,3}(\omega) = \frac{\eta_1(\omega)}{2} \mp \sqrt{\frac{\vartheta_1^2}{4} - 1} \; \cdot \sum_{n=0}^\infty \frac{(-1)^n \, (2n)! \, r_1^n(\omega)}{(1-2n) \, n! \, 4^n \, (\theta_1+2)^n} \cdot \sum_{n=0}^\infty \frac{(-1)^n \, (2n)! \, r_1^n(\omega)}{(1-2n) \, n! \, 4^n \, (\theta_1-2)^n},
\end{equation}
where $\eta_1(\omega)$ and $r_1(\omega)$ are taken from Eqs.~(\ref{eta1series}) and (\ref{4.10}), respectively.

In the same way we obtain the Taylor series expansion for other characteristic roots:
\begin{equation}
    \lambda_{2,4}(\omega) = \frac{\eta_2(\omega)}{2} \mp \sqrt{\frac{\vartheta_2^2}{4} - 1} \; \cdot \sum_{n=0}^\infty \frac{(-1)^n \, (2n)! \, r_2^n(\omega)}{(1-2n) \, n! \, 4^n \, (\theta_2+2)^n} \cdot \sum_{n=0}^\infty \frac{(-1)^n \, (2n)! \, r_2^n(\omega)}{(1-2n) \, n! \, 4^n \, (\theta_2-2)^n},
\end{equation}
where functions $\eta_2(\omega)$ and $r_2(\omega)$ are taken from Eqs.~(\ref{eta2series}) and (\ref{4.10}), respectively.

\section{{Legendre Polynomial Calculation}}
\label{app.legendre}
In our numerical experiments (in Subsect.~\ref{subse.TBC} and \ref{app.eta12}) we use coefficients of Taylor expansions of algebraic functions and their rational approximations (see \cite{nuttall1984asymptotics}). Convergence of such series depends on the location of functions' singular points on the complex plane.

We need to calculate Legendre polynomials $\mathrm{P}_n\left(\varepsilon\right)$ of order $n$ in Eqs.~\eqref{eta1series}, \eqref{eta2series}. It is necessary to assess the speed of decrease of $\mathrm{P}_n(\varepsilon)$ with increasing order $n$. Asymptotic formula for $\varepsilon = \cos \alpha$ as $n\to\infty$ (see, e.g. \cite{OLVER1974, abramowitz1972handbook})
\begin{equation}
    \label{app_eq:Legendre_asymp}
    \mathrm{P}_n(\cos\alpha)=\sqrt{\frac{2}{\pi n \sin\alpha}}\sin\left(n\alpha+\alpha/2 + \pi/4\right) + \mathbf{O}\left(n^{-3/2}\right)
\end{equation}
holds for $|\varepsilon| < 1$. In our case, we need to make sure that the absolute value of Legendre polynomials' argument is less than one: $\left|\frac{\mu^2}{\mu^2-2\nu}\right| < 1$, which is true if $\mu^2 < \nu$. By definition, it leads to
\begin{equation*}
    \tau > \sqrt{\frac{R^2\,\rho}{E}}.
\end{equation*}
Note that in our numerical experiment (see Subsect.~\ref{subse.rod}) it corresponds to $\tau > \approx 1.9346 \, \cdot 10^{-6} \, \mathrm{s}$. In some other examples the inequality $|\varepsilon| < 1$ coincides with the Courant stability condition of finite-difference scheme, see~\cite{gordin1979upper, gordin1987mathematicalb, gordin2000mathematical}.

From Eq.~\eqref{app_eq:Legendre_asymp} we see that the decrease of Legendre polynomial values is not fast. With increasing $n$ the calculation of Legendre polynomials becomes a numerically difficult task. Standard recurrent formula at point $\varepsilon$
\begin{equation*}
    (n+1) \mathrm{P}_{n+1}(\varepsilon) = (2n+1)x\mathrm{P}_n(\varepsilon) - n\mathrm{P}_{n-1}(\varepsilon)
\end{equation*}
may accumulate a big error for high $n$. We use the approach proposed by S.L.Belousov (see \cite{BELOUSOV1962ru, BELOUSOV19623}):
\begin{equation}
\label{app_eq_Belousov_P}
\begin{aligned}
    \mathrm{P}_n(\cos \alpha) = \sqrt{\frac{2n+1}{2}} \, \frac{1\cdot3\cdot5\ldots(2n-1)}{2^{n-1}\cdot n!} &\left[ \cos n\alpha + \frac{1}{1} \frac{n}{2n-1} \cos(n-2)\alpha +  \right. \\
    & + \frac{1\cdot3}{1\cdot2} \frac{n(n-1)}{(2n-1)(2n-3)} \cos(n-4)\alpha + \\
    & + \frac{1\cdot3\cdot5}{1\cdot2\cdot3} \frac{n(n-1)(n-2)}{(2n-1)(2n-3)(2n-5)} \cos(n-6)\alpha + \ldots \\
    \ldots & + 
    \left\{\begin{array}{rl}
        \frac{(n-2)!!}{\lfloor n/2 \rfloor!} \frac{\prod_{k=1}^{\lfloor n/2 \rfloor} (n-k+1)}{\prod_{k=1}^{\lfloor n/2 \rfloor}(2n - 2k + 1)} \cos \alpha & \text{if } n \text{ is odd,} \\
        \frac{1}{2} \cdot \frac{(n-2)!!}{(n/2)!} \frac{\prod_{k=1}^{n/2} (n-k+1)}{\prod_{k=1}^{n/2}(2n - 2k + 1)} \cos 0\alpha & \text{if } n \text{ is even}
    \end{array}\right].
\end{aligned}
\end{equation}

Here brackets $\lfloor \cdot \rfloor$ denote rounding down (floor function). The series in Eq.~\eqref{app_eq_Belousov_P} conclude with a term containing $\cos \alpha$ with $n$ being odd. If $n$ is even, the series is concluded with a term containing $\cos 0 \alpha$. The last coefficient (before $\cos 0 \alpha$) is additionally multiplied by $1/2$.

\section{Stability of the Crank -- Nicolson approximation of the rod transverse vibrations equation}
\label{app_CN_stability}
Here we investigate the stability of Crank -- Nicolson approximation of rod transverse vibrations equation on infinite domain, i.e. without boundary conditions. After the Fourier transform of Eq.\eqref{2.2}
\begin{equation*}
    \sigma \left(u_{m+2}^{n+1} + u_{m-2}^{n+1} + u_{m+2}^{n-1} + u_{m-2}^{n-1}\right) + \beta \left(u_{m+1}^{n+1} + u_{m-1}^{n+1} + u_{m+1}^{n-1} + u_{m-1}^{n-1} \right) + \alpha \left(u_m^{n+1} + u_m^{n-1}\right) + \gamma \left(u_{m-1}^{n} + u_{m+1}^{n}\right) + \delta u_m^n = 0
\end{equation*}
we get the ordinary finite-difference equation
\begin{equation}
\label{app_CNsymbol}
    \sigma \left[ 2\cos(2\xi h) \left(w^{n+1} + w^{n-1}\right) \right] + \beta \left[ 2\cos(\xi h) \left(w^{n+1} + w^{n-1}\right) \right] + \alpha \left(w^{n+1} + w^{n-1}\right) + \left[2\gamma\cos(\xi h) + \delta \right] w^n = 0,
\end{equation}
where $\xi$ is dual to spatial variable, $w$ is the solutions's image.

The equation for the next time step has the matrix form
\begin{equation}
\label{app_matHard}
    \begin{pmatrix}
        w^{n+1} \\
        w^n
    \end{pmatrix} = 
    \begin{pmatrix}
        -\frac{2\gamma \cos(\xi h) + \delta}{2\sigma \cos(2\xi h) + 2\beta \cos(\xi h) + \alpha} & -1 \\
        1 & 0
    \end{pmatrix}
    \begin{pmatrix}
        w^n \\
        w^{n-1}
    \end{pmatrix}.
\end{equation}

Simplifying Eq.~\eqref{app_matHard} we get
\begin{equation}
\label{app_mat}
    \begin{pmatrix}
        w^{n+1} \\
        w^n
    \end{pmatrix} = 
    \begin{pmatrix}
        -\frac{4\mu \cos(\xi h) -2-4\mu}{2\nu\cos^2(\xi h) - 2(2\nu+\mu) \cos(\xi h) + 1 + 2\nu + 2\mu} & -1 \\
        1 & 0
    \end{pmatrix}
    \begin{pmatrix}
        w^n \\
        w^{n-1}
    \end{pmatrix}.
\end{equation}

For scheme \eqref{2.2} to be stable, it is necessary for the largest absolute value of the eigenvalue of matrix in Eq.~\eqref{app_mat} to be less or equal to one for any real $\xi$:
\begin{equation}
    \max_{i = 1, 2} | \mathrm{eig}_i(\xi) | \leq 1, \quad \forall \xi \in \mathbb{R}.
\end{equation}

According to Vieta's theorem for characteristic equation for matrix in Eq.~\eqref{app_mat}, the product of both eigenvalues is equal to one. Therefore, the stability is obtained when the discriminant of the corresponding characteristic equation is non-positive:
\begin{equation}
\label{app_stab_cond}
    \left(\frac{2\mu \cos(\xi h) - 1 - 2\mu}{2\nu\cos^2(\xi h) - 2(2\nu+\mu) \cos(\xi h) + 1 + 2\nu + 2\mu} \right)^2 - 1 \leq 0.
\end{equation}

The numerator in the {parentheses} in Eq.~\eqref{app_stab_cond} is negative for any $\xi$:
\begin{equation}
    \label{app_numer}
    2\mu(\cos(\xi h) - 1) - 1 < 0.
\end{equation}

The denominator in {parentheses} in Eq.~\eqref{app_stab_cond} can be rewritten as
\begin{equation}
    2\nu \left(1 - \cos(\xi h) \right)^2 + 2\mu\left(1 - \cos(\xi h) \right) + 1 > 0 \; \forall \xi \in \mathbb{R}.
\end{equation}

Therefore, the stability condition Eq.~\eqref{app_stab_cond} becomes
\begin{equation*}
    -2\mu(\cos(\xi h) - 1) + 1 \leq  2\nu \left(1 - \cos(\xi h) \right)^2 + 2\mu\left(1 - \cos(\xi h) \right) + 1,
\end{equation*}
or simply
\begin{equation}
\label{app_stab}
    2\nu\left(1-\cos(\xi h)\right)^2 \geq 0,
\end{equation}
which is true for any real $\xi$. Hence, the approximation is absolutely stable.

\section{{Initial Data Construction}}
\label{app_init}
Let us decompose solution $u(t,x)$ into a Taylor series with respect to time in the vicinity of $t=0$:
\begin{equation}
\label{6.1}   
u(\tau,x)=U_0(x)+U_1(x)\tau+U_2(x)\frac{\tau^2}{2!}+\ldots ,
\end{equation}
where $U_k(x)=\frac{\partial^k u}{\partial t^k}(0,x)$. The functions $U_0,\,U_1$ compose initial conditions for Eq.~(\ref{1.1}). To determine the left hand side in Eq.~(\ref{6.1}) with an error $O\left(\tau^3\right)$, we need to compute the function $U_2$.

We differentiate Eq.~(\ref{6.1})

\begin{equation*}
\begin{aligned}
    \frac{\partial^2 u}{\partial t^2}(\tau, x) &= \sum_{k=2}^\infty U_k (x) \frac{\tau^{k-2}}{(k-2)!}, \\
    \frac{\partial^4 u}{\partial t^2 \, \partial x^2}(\tau, x) &= \sum_{k=2}^\infty U_k^{''} (x) \frac{\tau^{k-2}}{(k-2)!}, \\
    \frac{\partial^4 u}{\partial x^4}(\tau, x) &= \sum_{k=2}^\infty U_k^{[iv]} (x) \frac{\tau^k}{k!}, \\
\end{aligned}
\end{equation*}
and substitute the series into Eq.~(\ref{2.1}):
\begin{equation*}
    \sum\limits_{k=2}^\infty U_k(x)\frac{\tau^{k-2}}{(k-2)!} - R^2 \sum_{k=2}^\infty U_k^{''}(x) \frac{\tau^{k-2}}{(k-2)!} + \frac{E R^2}{\rho} \sum_{k=0}^\infty U_k^{[iv]}(x) \frac{\tau^k}{k!} = 0.
\end{equation*}

Collecting similar terms at the  zeroth degree of $\tau$ we obtain the following linear ordinary differential equation for $U_2$ with constant coefficients:
\begin{equation}
\label{6.2}   
\left[D\frac{\mathrm{d}^2}{\mathrm{d}x^2}-1\right]U_2(x)=CU^{[iv]}_0(x),
\end{equation}
where $D = R^2$ and $C = E\,R^2\,\rho^{-1}$.

We assume here that both initial function $U_0^{[iv]}(x)$ and auxiliary function $U_2(x)$ are rapidly decreasing at infinity. In this case, the principal term of the asymptotic at infinity of the {non-increasing} solution of the non-homogeneous differential equation, can be determined from the homogeneous one:
\begin{equation}
\label{6.3}   
U_2(x)\sim\exp\left(\frac{-x}{\sqrt{D}}\right)\;\mbox{as}\:x\to +\infty,
\quad U_2(x)\sim\exp\left(\frac{x}{\sqrt{D}}\right)\;\mbox{as}\:x\to -\infty.
\end{equation}

We approximate Eq.~(\ref{6.2}) by a compact finite-difference scheme (for details see, e.g. \cite{lele1992compact, gordin2010mathematics}):
\begin{equation}
\label{6.4}
\begin{aligned}
    a U_2(x_{j-1}) + U_2(x_j) + a U_2(x_{j+1}) &= p U_0(x_{j-2}) + q U_0(x_{j-1}) + r U_0(x_{j}) + q U_0(x_{j+1}) + p U_0(x_{j+2}),\\
    j &\in [-J+1,\,J-1].
\end{aligned}
\end{equation}

To determine the unknown four values $a$, $p$, $q$, $r$ we substitute into Eq.~(\ref{6.4}) four pairs of the test functions $\langle U_{2, k},\,U_{0, k} \rangle \equiv \langle u_k,\,f_k\rangle$, $k=0,\,1,\,2,\,3$, which satisfy Eq.~(\ref{6.2}) and are listed in Table~\ref{tab:1}.

\begin{table}[h]
    \centering
    \begin{tabular}{cccl}
    \toprule
    \textbf{No} & $u_k$ & $f_k$ & Algebraic equation for the coefficients of Eq.~(\ref{6.4})\\
    \midrule
      0 & 0 & 1 & $2p+2q+r=0$\\
      1 & 0 & $(x-x_j)^2$ & $8p+2q=0$\\
     2 & $-24 C$ &  $(x-x_j)^4$ & $-24C(2a+1)=(32p+2q)h^4$\\
     3 & $-360C((x-x_j)^2+2D)$ & $(x-x_j)^6$ & $-720DC(2a+1)-720Cah^2=2p(2h)^6+2qh^6$ \\
    \bottomrule
    \end{tabular}
    \caption{The pairs of test functions and the corresponding linear algebraic equations for coefficients of compact finite-difference scheme (\ref{6.4}).}
    \label{tab:1}
\end{table}

Differential Eq.~(\ref{6.2}) and finite-difference Eq.~(\ref{6.4}) both include only even order operators with constant coefficients. Therefore, the following odd 
test functions $\langle u,\,f\rangle$ like $\langle 0,\,(x-x_j) \rangle$ or $\langle (x-x_j)^{2m+1}$, $D\cdot (2m+1)\cdot 2m\cdot (x-x_j)^{2m-1}- (x-x_j)^{2m+1}\rangle$ satisfy Eq.~(\ref{6.4}) for any coefficients. That is why only even test functions are used for the approximation.

Solution of the system of these four linear algebraic equations is
\begin{equation*}
    a=\frac{h^2-6D}{12D+4h^2},\:
    p=\frac{-3C}{2h^2(3D+h^2)},\:
    q=\frac{6C}{h^2(3D+h^2)},\:
    r=\frac{-9C}{h^2(3D+h^2)}.
\end{equation*}

Three-diagonal system of $2J-1$ linear algebraic Eqs.~(\ref{6.4}) and the finite-difference approximations of  formulae (\ref{6.3}): $U_2(x_{-J})=\exp(-h/\sqrt{D})\, U_2(x_{-J+1})$ and $U_2(x_{J})=\exp(-h/\sqrt{D}) \, U_2(x_{J-1})$ is closed and non-degenerate. It can be solved by the classical double-sweep method.

If it is necessary to increase the order of accuracy of the initial function $u(\tau,x)$, we should similarly define a function $U_3 (x)$ and substitute it into Eq.~(\ref{2.1}), etc.

\section{{Possible Modification of ADTBCs}}
\label{app.modif}
The function $u(t,\, x) \equiv \text{const}$ is a solution of differential Eq.~(\ref{2.1}) and the finite-difference Eq.~\eqref{2.2}. We can additionally require the {DTBCs and ADTBCs} to satisfy this solution. In other words, we introduce an additional linear condition on the polynomial coefficients in Syst.~\eqref{5.4} and Syst.~\eqref{5.5} for $k=1$ and for $k=2$  (their total sums must be equal to zero):
\begin{equation}
\label{7.2.1}
    P_k(1) + Q_k(1) + R_k(1) + S_k(1) = 0.
\end{equation}

As mentioned in Subs.~\ref{subse.rod}, an even number of unknown coefficients is required to construct the {ADTBCs}. To maintain the same approximation order of Syst.~\eqref{5.4} and Syst.~\eqref{5.5} with additional condition (\ref{7.2.1}), one extra coefficient is required and, thus, the number of these coefficients becomes odd.

\begin{figure}[ht]
    \centering
    \begin{subfigure}[b]{.33\linewidth}
        \includegraphics[width=\textwidth]{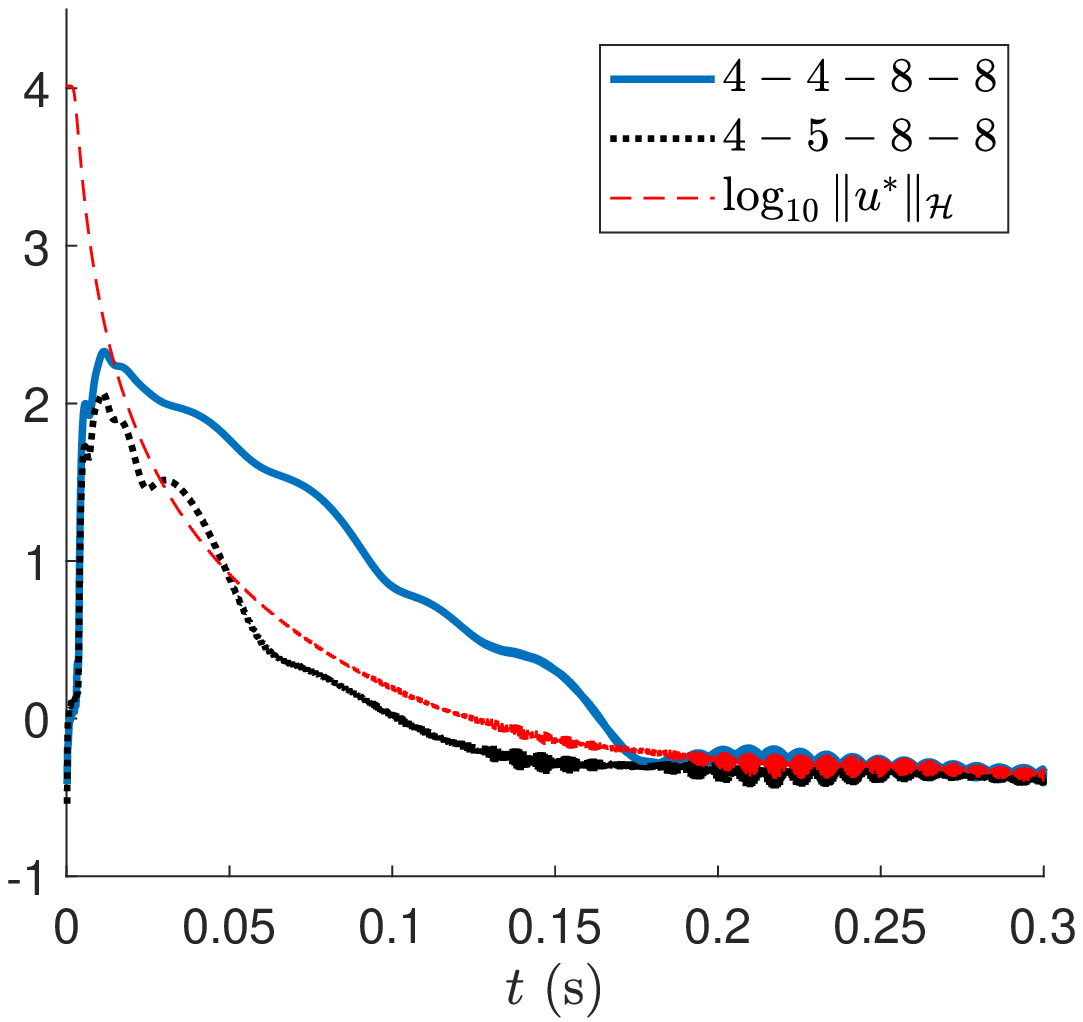}
        \caption{\label{fig:5a}}
    \end{subfigure}
    \begin{subfigure}[b]{.33\linewidth}
        \includegraphics[width=\textwidth]{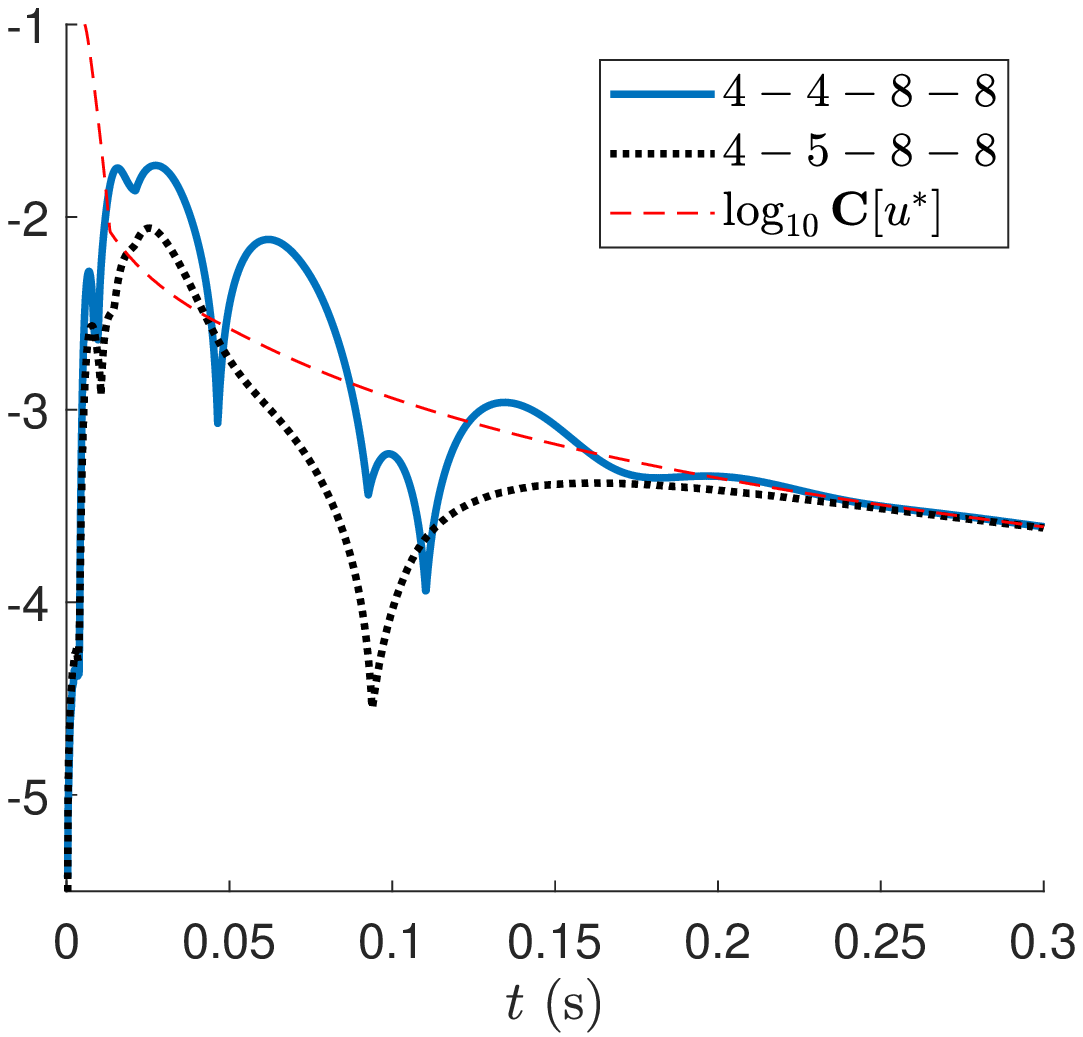}
        \caption{\label{fig:5b}}
    \end{subfigure}
    \begin{subfigure}[b]{.33\linewidth}
        \includegraphics[width=\textwidth]{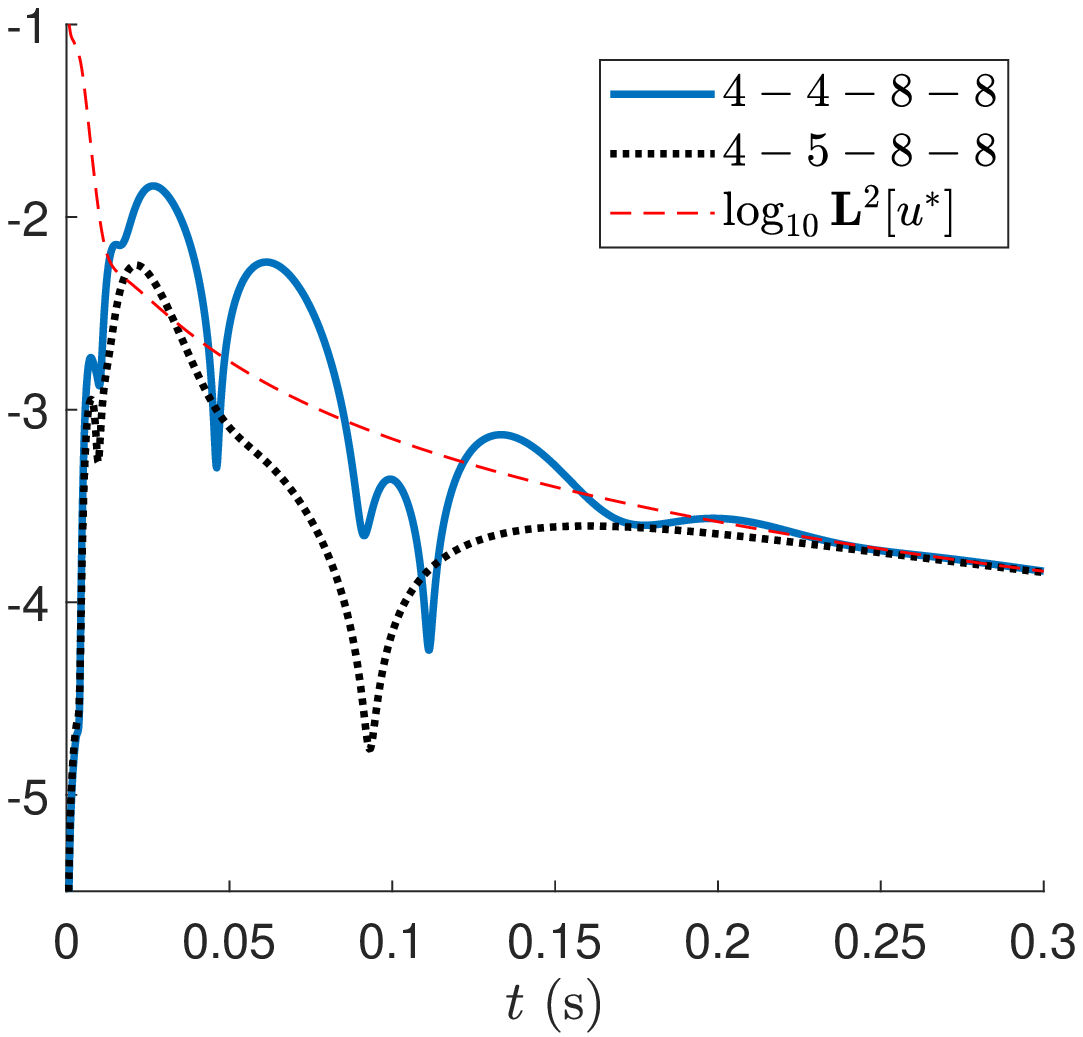}
        \caption{\label{fig:5c}}
    \end{subfigure}
    \caption{Solid and dotted lines represent {common logarithm of (a) {$\hat{\mathcal{H}}$}, (b) \textbf{C}-norm, (c) $\mathbf{L}^2$-norm of the difference between the reference solution $u^*$ and the solutions obtained with {ADTBCs}. The dotted line corresponds to the {ADTBCs} with additional requirement Eq.~(\ref{7.2.1}) taken into account in the polynomial coefficient calculation. The dash line stands for common logarithm of (a) {$\hat{\mathcal{H}}$}, (b) \textbf{C}-norm, (c) $\mathbf{L}^2$-norm of the reference solution $u^*$. Initial conditions are defined by formula~\eqref{init_cond_exp}}.}
    \label{fig:5}
\end{figure}

As a new example, we modify previous {ADTBCs} sets:
\begin{equation}
\label{7.2.2}
    \text{deg}\,P_k = 4, \; \text{deg}\,Q_k = 5, \; \text{deg}\,R_k = \text{deg}\,S_k = 8, \quad k=1,\,2.
\end{equation}

The coefficients of these sets that are derived from Syst.~\eqref{5.4} and Syst.~\eqref{5.5} with additional condition (\ref{7.2.1}) are presented in Table~\ref{tab:3}.

\begin{table}[t]
    \centering
    \begin{tabular}{l|rrrr | rrrr}
         & $P_1$ & $Q_1$ & $R_1$ & $S_1$ & $P_2$ & $Q_2$ & $R_2$ & $S_2$  \\
        \toprule
        1&1&	0&	-0.555979&	0.278657&   0&  1&  -0.925737&  0.301010 \\
        $\omega$&   -2.554692&  2.432054&   -1.468661&  0.329664&   -0.491692&  -0.454239&  0.247374&   -0.084630 \\
        $\omega^2$& 2.067232&   -1.792876&  0.313255&   0.091235&    0.249452&   0.758283&   -0.722172&  0.255585 \\
        $\omega^3$& -2.170815&  2.376252&   -0.936209&  0.136683&    -0.262835&	-0.272165&	0.412516&	-0.146529 \\
        $\omega^4$& 1.325085&   -1.388514&	0.900316&	-0.202468&   0.322920&	-0.303600&	0.219016&	-0.052243 \\
        $\omega^5$&         &   -0.519196&  0.545430&	-0.195266&           &	-0.155023&	0.171840&	-0.063135 \\
        $\omega^6$&         &	        &   0.009158&	-0.011746&           &	        &	0.004913&	-0.005389 \\
        $\omega^7$&         &	        &	-0.007229&	0.000655&            &	        &	-0.002735&	0.000079 \\
        $\omega^8$&         &	        &	-0.003216&	0.001191&            &	        &	-0.001341&	0.000474 \\
    \end{tabular}
\caption{The coefficients of the {ADTBCs} are obtained from Eqs.~(\ref{2zv}) and (\ref{3zv}), and additional condition (\ref{7.2.1}) for the sets of polynomial degrees $\langle P_k, \, Q_k, \, R_k, \, S_k \rangle = \langle 4, \, 5, \, 8, \, 8 \rangle$ at $k = 1,\,2$.}
\label{tab:3}
\end{table}

{
The introduced modification of DTBCs and ADTBCs construction can only be applied when the initial function $U_0 \equiv u(0, x)$ has zeroth integral over the segment:}
\begin{equation}
\label{eq:zeroth_int}
    \int_{-L/2}^{L/2} u(0, x) \, \mathrm{d}x = 0.
\end{equation}

\begin{figure}[ht]
    \centering
    \begin{subfigure}[b]{.33\linewidth}
        \includegraphics[width=\textwidth]{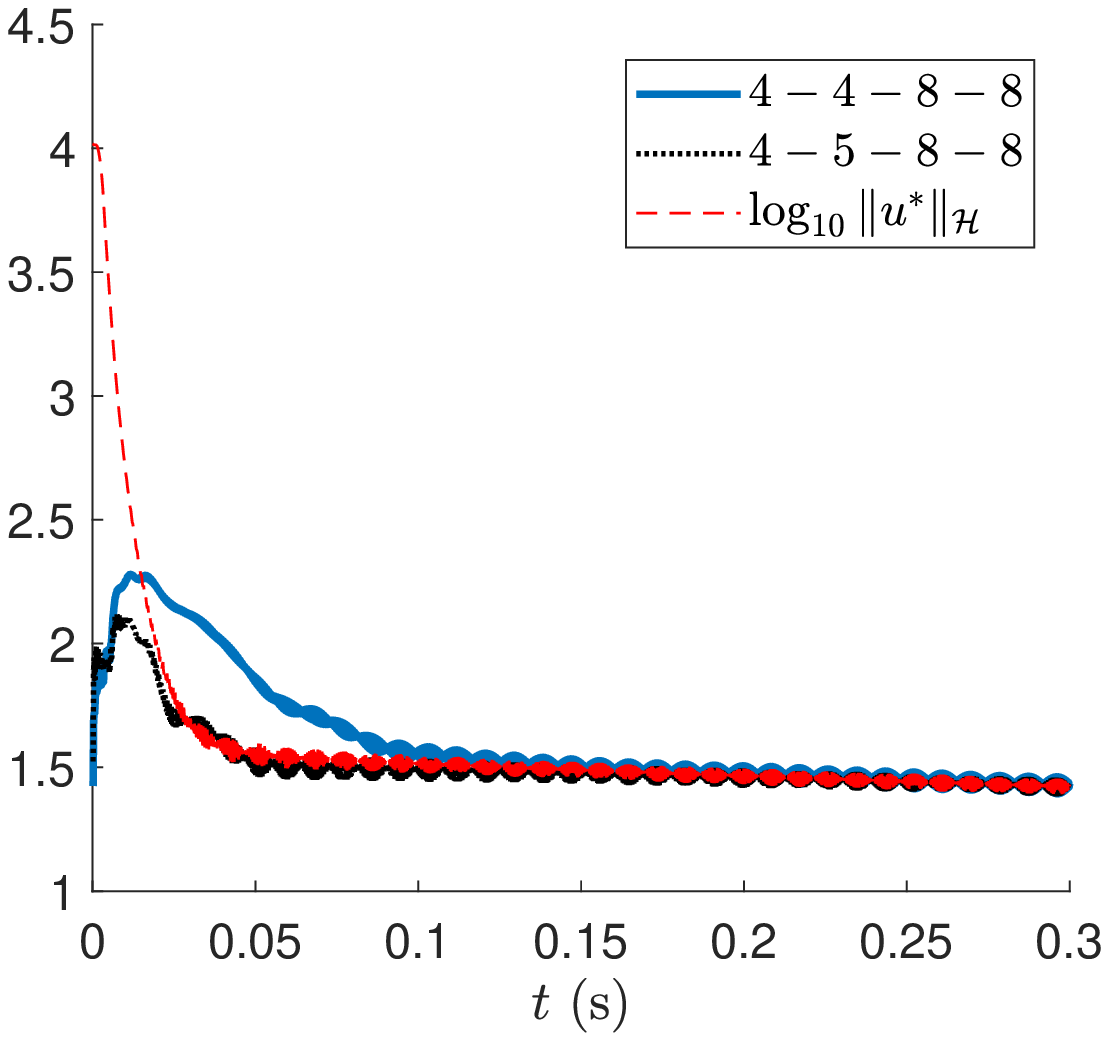}
        \caption{\label{fig:6a}}
    \end{subfigure}
    \begin{subfigure}[b]{.33\linewidth}
        \includegraphics[width=\textwidth]{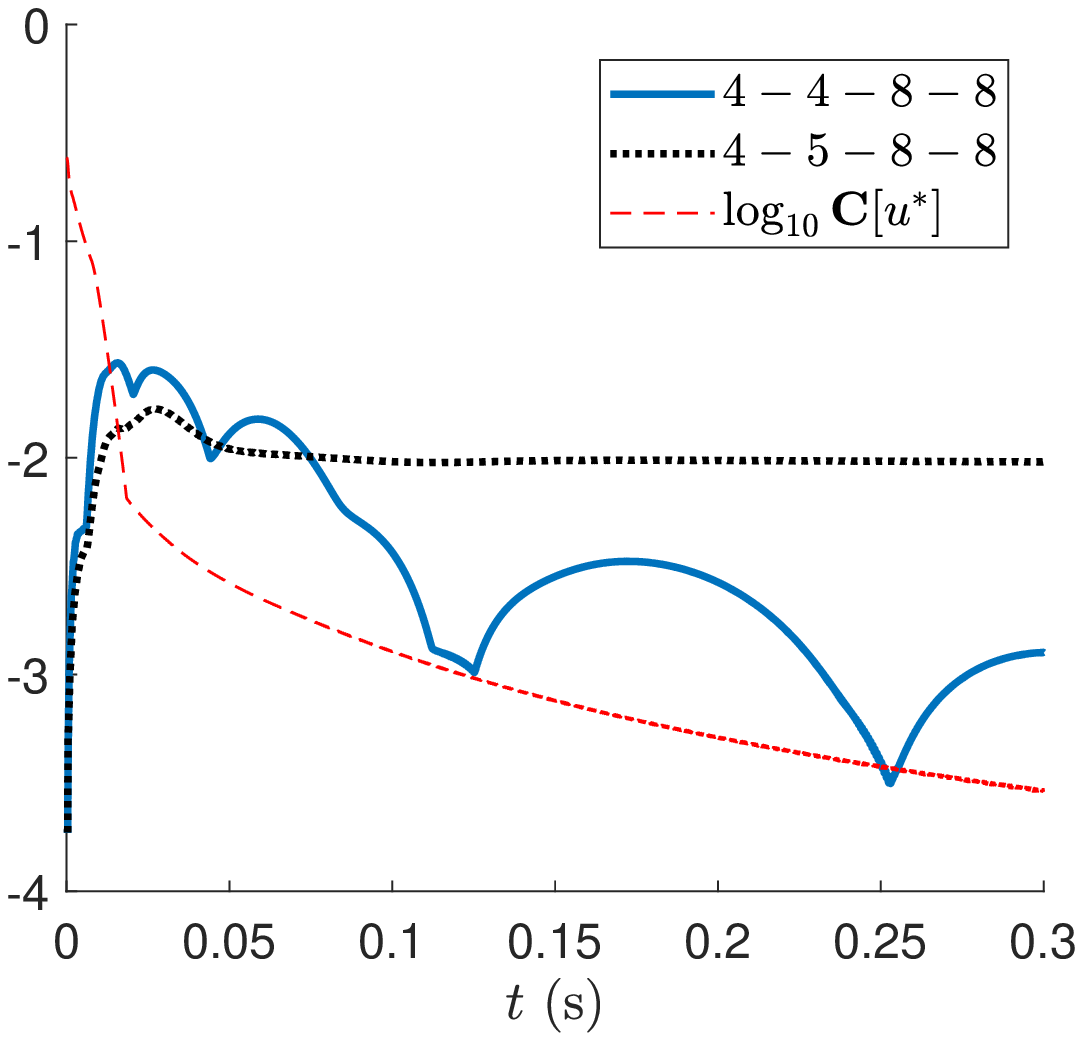}
        \caption{\label{fig:6b}}
    \end{subfigure}
    \begin{subfigure}[b]{.33\linewidth}
        \includegraphics[width=\textwidth]{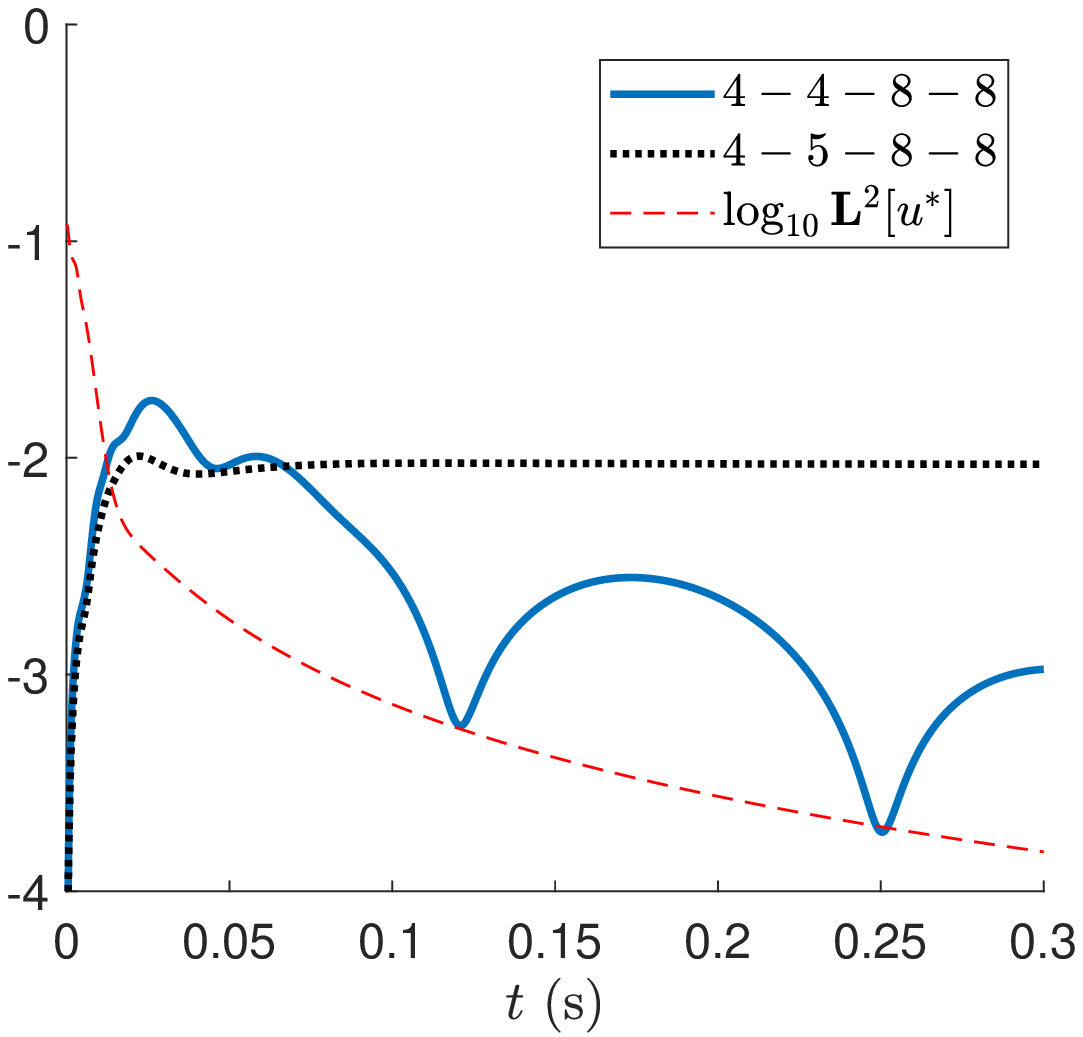}
        \caption{\label{fig:6c}}
    \end{subfigure}
    \caption{{Solid and dotted lines represent common logarithm of (a) {$\hat{\mathcal{H}}$}, (b) \textbf{C}-norm, (c) $\mathbf{L}^2$-norm of the difference between the reference solution $u^*$ and the solutions obtained with {ADTBCs}. The dotted line corresponds to the {ADTBCs} with additional requirement Eq.~(\ref{7.2.1}) taken into account in the polynomial coefficient calculation. The dash line is common logarithm of (a) {$\hat{\mathcal{H}}$}, (b) \textbf{C}-norm, (c) $\mathbf{L}^2$-norm of the reference solution $u^*$. Initial conditions are defined by formula~\eqref{eq:u0-nonsym}.}}
    \label{fig:6}
\end{figure}

{
When the condition~\eqref{eq:zeroth_int} does not hold, the resulting {ADTBCs} do not provide the transparency property. In Fig.~\ref{fig:6}, we present the errors for sets $\langle 4,4,8,8\rangle$ and $\langle 4,5,8,8\rangle$, and shifted initial conditions:}
\begin{equation}
\label{eq:u0-nonsym}
    u(0, x) = \frac{(x - 0.1)}{{\sqrt{\pi\cdot 0.02}}} \cdot \exp\left(-\frac{(x - 0.1)^2}{0.02}\right), \quad \partial_t u(0, x) = 0, \quad x \in \left[-L/2,\, L/2\right].
\end{equation}

Approximate solution's error under the {ADTBCs} (\ref{7.2.2}) with additional requirement (\ref{7.2.1}) is mostly smaller (see Fig.~\ref{fig:5}). Therefore, the introduction of this condition on coefficients results in a decrease of the error for any time moment $t${, provided the initial conditions satisfy Eq.~\eqref{eq:zeroth_int}.}

\section{Dynamics of obtained solutions for different {ADTBCs}}
\label{app.anim}
Here we present animations of the obtained solutions $u$ and the reference solution $u^*$ (left). We also provide the common logarithm of $\mathbf{C}$-norm of absolute difference between the two solutions (right).

Sets of polynomial degrees are denoted as $\mathbf{deg} P_k - \mathbf{deg} Q_k - \mathbf{deg} R_k - \mathbf{deg} S_k$ (being equal for $k = 1$ and $k = 2$) in the title of animations.

All results are obtained using the same parameters and initial conditions as in Subsect.~\ref{subse.rod}.





\newpage
\section*{References}
\bibliographystyle{elsarticle-num-names}






\end{document}